\documentclass[10pt,twoside]{article}
\usepackage{amsfonts} 

\font\caps=cmb10                    
\font\Caps=cmbx10 scaled \magstep1   

\def\firstpage{1}\def\lastpage{1000}

\makeatletter
\@specialpagefalse\headheight=8.5pt
\def\DocMath{}
\renewcommand{\@evenhead}{%
    \ifnum\thepage>\lastpage\rlap{\thepage}\hfill%
    \else\rlap{\thepage}\slshape\leftmark\hfill\caps\SAuthor\hfill\fi}%
\renewcommand{\@oddhead}{%
    \ifnum\thepage=\firstpage{\DocMath\hfill\llap{\thepage}}%
    \else{\slshape\rightmark}\hfill\caps\STitle\hfill\llap{\thepage}\fi}%
\makeatother

\def\TSkip{\bigskip}
\newbox\TheTitle{\obeylines\gdef\GetTitle #1
\ShortTitle  #2
\SubTitle    #3
\Author      #4
\ShortAuthor #5
\EndTitle
{\setbox\TheTitle=\vbox{\baselineskip=20pt\let\par=\cr\obeylines%
\halign{\centerline{\Caps##}\cr\noalign{\medskip}\cr#1\cr}}%
	\copy\TheTitle\TSkip\TSkip%
\def\next{#2}\ifx\next\empty\gdef\STitle{#1}\else\gdef\STitle{#2}\fi%
\def\next{#3}\ifx\next\empty%
    \else\setbox\TheTitle=\vbox{\baselineskip=20pt\let\par=\cr\obeylines%
    \halign{\centerline{\caps##} #3\cr}}\copy\TheTitle\TSkip\TSkip\fi%
\centerline{\caps #4}\TSkip\TSkip%
\def\next{#5}\ifx\next\empty\gdef\SAuthor{#4}\else\gdef\SAuthor{#5}\fi%

\TSkip\TSkip%
\catcode'015=5}}\def\Title{\obeylines\GetTitle}
\def\Abstract{\begingroup\narrower
    \parskip=\medskipamount\parindent=0pt{\caps Abstract. }}
\def\EndAbstract{\par\endgroup\TSkip\TSkip}
\newbox\TheAdd\def\Addresses{\vfill\copy\TheAdd\vfill
    \ifodd\number\lastpage\vfill\eject\phantom{.}\vfill\eject\fi}
{\obeylines\gdef\GetAddress #1
\Address #2 
\Address #3
\Address #4
\EndAddress
{\def\xs{6truecm}
\setbox0=\vtop{{\obeylines\hsize=\xs#1}}\def\next{#2}
\ifx\next\empty 
     \setbox\TheAdd=\hbox to\hsize{\hfill\copy0\hfill}
\else\setbox1=\vtop{{\obeylines\hsize=\xs#2}}\def\next{#3}
\ifx\next\empty 
     \setbox\TheAdd=\hbox to\hsize{\hfill\copy0\hfill\copy1\hfill}
\else\setbox2=\vtop{{\obeylines\hsize=\xs#3}}\def\next{#4}
\ifx\next\empty\ 
     \setbox\TheAdd=\vtop{\hbox to\hsize{\hfill\copy0\hfill\copy1\hfill}
	        \vskip20pt\hbox to\hsize{\hfill\copy2\hfill}}
\else\setbox3=\vtop{{\obeylines\hsize=\xs#4}}
     \setbox\TheAdd=\vtop{\hbox to\hsize{\hfill\copy0\hfill\copy1\hfill}
	        \vskip20pt\hbox to\hsize{\hfill\copy2\hfill\copy3\hfill}}
\fi\fi\fi\catcode'015=5}}\gdef\Address{\obeylines\GetAddress}

\begin{document} 

\Title Flat connections, Higgs operators, and Einstein metrics 
on compact Hermitian manifolds
\ShortTitle Flat connections, Higgs operators, Einstein metrics
\SubTitle   
\Author M. L\"ubke
\ShortAuthor
\EndTitle

\Abstract 
A flat complex vector bundle $(E,D)$ on a compact Riemannian
manifold $(X,g)$ is stable (resp. polystable) in the sense of Corlette [C] if
it has no $D$-invariant subbundle (resp. if it is the $D$-invariant direct sum
of stable subbundles). It has been shown in [C] that the polystability of
$(E,D)$ in this sense is equivalent to the existence of a so-called harmonic
metric in $E$. In this paper we consider flat complex vector bundles on compact
Hermitian manifolds $(X,g)$. We propose new notions of $g$-(poly-)stability
of such bundles, and of $g$-Einstein metrics in them; these notions
coincide with (poly-)stability and harmonicity in the sense of Corlette if
$g$ is a K\"ahler metric, but are different in general. Our main result is
that the $g$-polystability in our sense is equivalent to the existence
of a $g$-Hermitian-Einstein metric. Our notion of a $g$-Einstein metric in a flat
bundle is motivated by a correspondence between flat bundles and Higgs bundles
over compact surfaces, analogous to the correspondence in the case of K\"ahler
manifolds [S1], [S2], [S3]. 

\EndAbstract

\Address M. L\"ubke
Mathematical Institute
Leiden University
PO Box 9512
NL 2300 RA LEIDEN
\Address
\Address
\Address
\EndAddress
%
%

\newcommand{\dpr}{^{\prime\prime}}
\newcommand{\Vol}{\mathrm{Vol}}
\newcommand{\rk}{\mathrm{rk}}
\newcommand{\eps}{\varepsilon}
\newcommand{\map}{\longrightarrow}
\newcommand{\mbb}{\mathbb}
\newcommand{\dbar}{\bar{\partial}}
\newcommand{\sett}[2]{\{\ #1\ |\ #2\ \}}
\newcommand{\End}{\mathrm{End}}
\newcommand{\Hom}{\mathrm{Hom}}
\newcommand{\half}{\frac{1}{2}}
\newcommand{\tr}{\mathrm{tr}}
\newcommand{\id}{\mathrm{id}}
\newcommand{\im}{\mathrm{im}}
\newcommand{\Pic}{\mathrm{Pic}}
\newcommand{\qed}{\hfill\vrule height6pt width6pt depth0pt \bigskip}
\newcommand{\qmod}[2]{{\hbox{}^{\displaystyle{#1}}}\!\big/\!\hbox{}_{\displaystyle{#2}}}
\newcommand{\mb}[1]{\hbox{$#1$}}

\newcommand{\ovrightarrow}[1]{\mathop{\vbox{\ialign{
				##\crcr
    ${\scriptstyle\hfil\;\;#1\;\;\hfil}$\crcr
    \noalign{\kern-1pt\nointerlineskip}
    \rightarrowfill\crcr}}\;}}

\newtheorem{thm}{Theorem}[section]
\newtheorem{prop}[thm]{Proposition}
\newtheorem{lem}[thm]{Lemma}
\newtheorem{cor}[thm]{Corollary}
\newtheorem{rem}[thm]{Remark}
\newtheorem{defn}[thm]{Definition}

\section{Introduction.}

Let $X$ be an $n$-dimensional compact complex manifold. If $X$ admits a
\underbar{K\"ahler} metric  $g$, then it is known by work of in particular
Simpson [S1],[S2],[S3] that there exists an  canonical identification of the
moduli space of polystable (or semisimple) flat bundles on $X$ with the moduli
space of $g$-polystable Higgs-bundles with vanishing Chern classes on $X$. 
This identification has been used in showing that certain groups are not
fundamental groups of compact K\"ahler manifolds. The construction uses the
existence of {\it canonical} metrics, called $g$-{\it harmonic} in the case of
flat bundles, and $g$-{\it Einstein} in the case of Higgs bundles. 

For flat bundles, the equivalence of semisimplicity and the existence
of a $g$-harmonic metric holds on compact \underbar{Riemannian} manifolds [C]. 
Furthermore, the equivalence of $g$-poly\-stability and the existence of a
$g$-Einstein metrics for Higgs bundles should generalize to the case of
\underbar{Hermitian} manifolds as in the case of holomorphic vector bundles,
using \underbar{Gauduchon} metrics. Nevertheless,
an identification as above cannot be expected for general compact Hermitian
manifolds, since it should imply restrictions on the fundamental group, but
every finitely presented group is the fundamental group of a 3-dimensional
compact complex manifold by a theorem of Taubes~[T].

In the case of compact complex \underbar{surfaces}, however, things are
different. We show that for an integrable Higgs
bundle $(E,d\dpr)$ with vanishing real Chern numbers and of $g$-degree 0 with
$g$-Einstein metric $h$ on a compact complex surface $X$ with Hermitian metric
$g$, there is an canonically associated flat connection $D$ in $E$, again of
$g$-degree 0, such that $h$ is what we call a $g$-\underbar{Einstein}
metric for $(E,D)$, and that the converse is also true. Furthermore, this
correspondence preserves isomorphism types and hence descends to a bijection
between moduli spaces. 

The notion of a $g$-Einstein metric in a flat bundle makes sense in higher
dimension, too, is equivalent to $g$-harmonicity in the case of a K\"ahler
metric, but different in general, and we show that the existence of such a
metric in a flat bundle $(E,D)$ is equivalent to the $g$-polystability of this
bundle in the sense that $E$ is the direct sum of $D$-invariant $g$-stable flat
subbundles. Here we call a flat bundle $(E,D)$ $g$-stable if every
$D$-invariant subbundle has $g$-\underbar{slope} larger(!) than the $g$-slope
of $(E,D)$. $g$-stability of a flat bundle is
equivalent to its stability (in the sense of Corlette) in the K\"ahler case,
but a weaker condition in general: A stable bundle is always $g$-stable,
but the tangent bundles of certain Inoue surfaces are examples of $g$-stable
bundles which are not stable. 

We expect that for a non-K\"ahler surface with Hermitian metric $g$, there is a
natural bijection between the moduli space of $g$-polystable Higgs bundles, with
vanishing Chern numbers and $g$-degree, and the moduli space of $g$-polystable
flat bundles with vanishing $g$-degree. In the last section we consider the
special case of \underbar{line} bundles on surfaces. Here the stability is
trivial, and the existence of Einstein metrics is easy to show, so we get indeed
the expected natural bijection between moduli spaces of line bundles of degree
0. We further show how this can be extended (in a non-natural way) to the moduli
spaces of line bundles of arbitrary degree; this extension argument works in
fact for bundles of arbitrary rank once the correspondence for degree 0 has been
established.

\vspace{.5cm}

{\bf Acknowledgments.} The author wishes to thank A. Teleman for the
suggestion to study the relation between flat connections and 
Higgs operators on Hermitian manifolds, and for several useful hints.  
Discussions with him and Ch. Okonek have been of considerable help in 
preparing this paper.\\
The author was supported in part by EC-HCM project AGE "Algebraic 
Geometry in Europe", contract no. ERBCHRXCT 940557.

\section{Preliminaries.}

Let $X$ be a compact $n$-dimensional complex manifold, and\ \ $E \map X$\ \ a
differentiable $\mbb{C}^r$-vector bundle on $X$. We fix the following

{\bf Notations:}\\
$A^p(X)$ (resp. $A^{p,q}(X)$) is the space of differentiable
$p$-forms (forms of type $(p,q)$) on $X$.\\
$A^p(E)$, $A^{p,q}(E)$ are the spaces of differential forms with values in $E$.\\
$\mathcal{A}(E)$ is the space of linear connections $D$ in $E$. For a connection\ \ 
$D \in \mathcal{A}(E)$\ \ we write\ \ $D = D^\prime + D\dpr\ ,$\ where 
 $D^\prime$ is of type (1,0) and $D\dpr$ of type (0,1). \\ 
$\mathcal{A}(E,h) \subset \mathcal{A}(E)$ is the subspace of $h$-unitary connections $d$ in $E$,
where $h$ is a Hermitian metric in $E$. We write\ \ $d = \partial + \dbar\ ,$\
where 
 $\partial$ is of type (1,0) and $\dbar$ of type (0,1).\\
$\mathcal{A}_f(E) := \sett{D \in \mathcal{A}(E)}{D^2 = 0}$\ \ is the subset of \underbar{flat}
connections. \\
$\bar{\mathcal{A}}(E)$ is the space of semiconnections $\dbar$ of type (0,1) in $E$
(i.e. $\dbar$ is the (0,1)-part of some\ \ $D \in \mathcal{A}(E)$\ ).\\
$\mathcal{H}(E) := \sett{\dbar \in \bar{\mathcal{A}}(E)}{\dbar^2 = 0}$\ \ is the subset of
\underbar{integrable} semiconnections or \underbar{holomorphic}
\underbar{structures} in $E$.\\
$\mathcal{A}^{\prime\prime}(E) := \bar{\mathcal{A}}(E)\oplus A^{1,0}(\End E) = \sett{d\dpr =
\dbar + \theta}{\dbar \in \bar{\mathcal{A}}(E)\ ,\ \theta \in A^{1,0}(\End E)}$ is the
space of \underbar{Higgs} \underbar{operators} in $E$. \\
$\mathcal{H}^{\prime\prime}(E) := \sett{d\dpr \in \mathcal{A}^{\prime\prime}(E)}{(d\dpr)^2 =
0}$ is the subset of \underbar{integrable} Higgs operators.\\
Often the same symbol is used for a connection, semiconnection, Higgs
operator etc. in $E$ and the induced operator in $\End E$.\\
Two connections\ \ $D_1,D_2 \in \mathcal{A}(E)$\ \ are \underbar{isomorphic},\
$D_1 \cong D_2\ ,$\ if there exists a differentiable automorphism $f$ of $E$
such that\ \ $f\circ D_1 = D_2\circ f\ ,$\ which is equivalent to\ \
$D(f) = 0\ ,$\ where $D$ is the connection in $\End E$ induced by $D_1$ and
$D_2$, i.e.\ \ $D(f) = D_2\circ f - f\circ D_1\ .$\ In the same way the
isomorphy of semiconnections resp. Higgs operators is defined. \\
If a Hermitian metric $h$ in $E$ is given, then a superscript ${}^*$ means adjoint
with respect to $h$. 

\vspace{.5cm}

For\ \ $D = D^\prime + D\dpr$\ \ there are unique
semiconnections $\delta_h^\prime$,$\delta_h\dpr$ of type (1,0), (0,1)
respectively such that $D^\prime + \delta_h\dpr$ and $\delta_h^\prime + D\dpr$
are $h$-unitary connections. Define\ \ $\delta_h := \delta_h^\prime +
\delta_h\dpr\ ;$\ then\ \ $d_h := \half(D + \delta_h)$\ \ is $h$-unitary, and\
\ $\Theta_h := D - d_h = \half(D - \delta_h)$\ \ is a $h$-selfadjoint 1-form with values in $\End E$. Let\ \ $d_h = \partial_h + \dbar_h$\ \ be the decomposition in the parts of type (1,0) and (0,1), and let $\theta_h$ be the $(1,0)$-part of $\Theta_h$; then it
holds
\[
D = d_h + \Theta_h = \partial_h + \dbar_h + \theta_h + \theta_h^*\ .
\]
The map 
\[
I_h : \mathcal{A}(E) \map \mathcal{A}\dpr(E)\ \ ,
\ \ I_h(D) := d\dpr_h := \dbar_h + \theta_h \in \mathcal{A}\dpr(E)
\]
is bijective; the inverse is given as follows. For\ \ $d\dpr = \dbar + \theta
\in \mathcal{A}\dpr(E)$\ \ let $\partial_h$ be the unique semiconnection of type $(1,0)$
such that the connection\ \ $d_h := \partial_h + \dbar$\ \ is $h$-unitary,
and define\ \ $\Theta := \theta + \theta^*\ .$\ Then
\[I_h^{-1}(d\dpr) = D_h := d_h + \Theta \in \mathcal{A}(E)\ .\]

\begin{rem}\label{rem1}
i) In general, if\ \ $D_1,D_2 \in \mathcal{A}(E)$\ \ are isomorphic, then\ \
$I_h(D_1)$ and $I_h(D_2)$ are \underbar{not} isomorphic, and
vice versa.\\
ii)\ $D_h = d_h + \theta +
\theta^*$\ \ is \underbar{not} $h$-unitary unless\ \ $\theta = 0\ ,$\ but
the connections\ \ $d_h - \theta + \theta^*$\ \ and
\ \ $d_h + \theta - \theta^*$\ \ are.\\
iii) Any metric $h^\prime$ in $E$ is of the form\ \ $h^\prime = f\cdot h\
,$\ i.e.\ \ $h^\prime(s,t) = h(f(s),t)\ ,$\ where $f$ is a $h$-selfadjoint
and positive definite. For a connection $D$ it is easy to show that the
operator $\delta_{h\cdot f}$ associated to $D$ and $f\cdot h$ is given by\ \
$\delta_{h\cdot f} = f^{-1}\circ\delta_h\circ f =
\delta_h + f^{-1}\circ\delta_h(f)\ ,$\ so it holds
\begin{eqnarray*}
d\dpr_{f\cdot h} &=& d\dpr_h + \half f^{-1}\circ\delta_h\dpr(f) - 
f^{-1}\circ\delta_h^\prime(f) \\
&=& d\dpr_h + \half f^{-1}\circ\dbar_h(f) - \half f^{-1}\circ\theta_h^*(f)  
- \half f^{-1}\circ\partial_h(f) + \half f^{-1}\circ\theta_h(f)\ .
\end{eqnarray*}
Conversely, for a given Higgs operator $d\dpr$ one verifies
\[
D_{f\cdot h} = D_h + f^{-1}\circ\partial_h(f) + f^{-1}\circ \theta(f)\ .
\] 
In particular, if $f$ is constant then the two maps $I_h$ and $I_{f\cdot h}$
coincide. 
 \end{rem}

\begin{defn} 
i) $G_h := (d\dpr_h)^2$\ \ is called the \underbar{pseudocurvature} of $D$ with
respect to $h$. \\
ii) $F_h := D_h^2$\ \ is called the \underbar{curvature} of $d\dpr$ with respect
to $h$. \end{defn}

\begin{rem}\label{rem2}
i) Obviously it holds:
$I_h(D)$ is an \underbar{integrable} Higgs operator if and only if\ \ $G_h 
= 0\ ,$\ and $I_h^{-1}(d\dpr)$ is a \underbar{flat} connection if and only if\ 
\ $F_h = 0\ .$\\
ii) For\ \ $i = 1,2\ ,$\ let $E_i$ be a differentiable complex vector
bundle on $X$ with Hermitian metric $h_i$ and connection $D_i$.
Let $h$ be the induced metric and $D$ the induced connection in $\Hom(E_1,E_2)$.  
Denote by $G_{i,h}$ resp. $G_h$ the pseudocurvature of
$D_i$ resp. $D$ with respect to $h_i$ resp. $h$. Then for\ \
$f \in A^0(\Hom(E_1,E_2))$\ \ it holds\ \ $G_h(f) =
G_{2,h}\circ f - f\circ G_{1,h}\ .$\\
Similarly, the curvature
$F_h$ of the Higgs operator induced in $\Hom(E_1,E_2)$ by Higgs
ope\-rators $d\dpr_i$ in the $E_i$ is given by\ \ $F_h(f) =
F_{2,h}\circ f - f\circ F_{1,h}\ .$\\
iii) If $D$ is a connection, then $D^2$ is the curvature of $d\dpr_h$
with respect to $h$, and if $d\dpr$ is a Higgs operator, then $(d\dpr)^2$ is
the pseudocurvature of $D_h$ with respect to $h$. This trivially follows from
the bijectivity of $I_h$. \end{rem}

\begin{lem}\label{lem1}
i) For\ \ $D \in \mathcal{A}(E)$\ \ let\ \ $D = d_h + \Theta_h = \partial_h +
\dbar_h + \theta_h + \theta_h^*$\ \ be the decomposition induced by $h$ as
above.
If $D$ is flat, then it holds\ \ $\delta_h^2 = 0$\ ,\ $d_h(\Theta_h)
= 0\ ,$\ i.e.\ \ $\partial_h(\theta_h) = \dbar_h(\theta_h^*) =
\partial_h(\theta_h^*) + \dbar_h(\theta_h) = 0\ ,$\ and furthermore\ \
$d_h^2 = -\Theta_h\wedge\Theta_h\ .$\\
ii) For\ \ $d\dpr = \dbar + \theta \in \mathcal{A}\dpr(E)$\ \ let $\partial_h$, $d_h$
and $D_h$ be as above, and write\ \ $d_h^\prime := \partial_h + \theta^*\ .$\\ 
If $d\dpr$ is integrable, then it holds\ \ $(d_h^\prime)^2 = 0\ 
,$\ i.e.\ \ $\partial_h^2 = \partial_h(\theta^*) = \theta^*\wedge\theta^* =
0\  ,$\ $d_h^2 = [\partial_h,\dbar]\ ,$\ and hence\ \ $F_h = d_h^2 +
[\theta,\theta^*] + \partial_h(\theta) + \dbar(\theta^*)\ .$ \end{lem}

{\bf Proof:} i) For\ \ \mb{D = D^\prime + D\dpr \in \mathcal{A}_f(E)}\ \ it holds
\begin{eqnarray*}
0 &=& \partial\partial h(s,t) = h((D^\prime)^2(s),t) - h(D^\prime(s),\delta\dpr_h(t)) + h(D^\prime(s),\delta\dpr_h(t)) + h(s,(\delta\dpr_h)^2(t))\\ &=& h(s,(\delta\dpr_h)^2(t))
\end{eqnarray*}
for all\ \ \mb{s,t \in A^0(E)\ ,} i.e.\ \ \mb{(\delta\dpr_h)^2 = 0\ .} Similarly one sees\ \ \mb{(\delta_h^\prime)^2 = 0 = \delta^\prime_h\delta\dpr_h + \delta\dpr_h\delta^\prime_h\ ,} yielding\ \ \mb{\delta_h^2 = 0\ .} We conclude
\[
d_h(\Theta_h) = \frac{1}{4}[D + \delta_h , D - \delta_h] = 0\ ,
\]
and
\[
0 = D^2 = (d_h + \Theta_h)^2 = d_h^2 + d_h(\Theta_h) + \Theta_h\wedge\Theta_h = d_h^2 + \Theta_h\wedge\Theta_h\ .
\]

ii) For\ \ \mb{d\dpr = \dbar + \theta \in \mathcal{H}\dpr(E)}\ \ and\ \ \mb{d_h = \partial_h + \dbar}\ \ it is well known that\ \ \mb{\partial_h^2 = 0\ ,} and hence\ \ \mb{d_h^2 = [\partial_h,\dbar]\ .} Furthermore, for all\ \ \mb{s,t \in A^0(E)}\ \ it holds
\begin{eqnarray*}
h(\partial_h(\theta^*)(s),t) &=& h(\partial_h\circ\theta^*(s),t) + h(\theta^*\circ\partial_h(s),t)\\
&=& \partial h(\theta^*(s),t) + h(\theta^*(s),\dbar(t)) - h(\partial_h(s),\theta(t))\\
&=& \partial h(s,\theta(t)) + h(s,\theta\circ\dbar(t)) - h(\partial_h(s),\theta(t))\\
&=& h(\partial_h(s),\theta(t)) + h(s,\dbar\circ\theta(t)) 
+ h(s,\theta\circ\dbar(t)) - h(\partial_h(s),\theta(t))\\
&=& h(s,\dbar(\theta)(t)) = 0\ ,
\end{eqnarray*}
and
\[
h(\theta^*\wedge\theta^*(s),t) = - h(s,\theta\wedge\theta(t)) = 0\ ;
\]
this shows\ \ \mb{\partial_h(\theta^*) = 0 = \theta^*\wedge\theta^*\ .}
\qed

Now let $g$ be a Hermitian metric in $X$, and denote by $\omega_g$ the 
associated $(1,1)$-form on $X$, by $\Lambda_g$ the 
contraction by $\omega_g$, and by $*_g$ the associated 
Hodge-$*$-operator. 

Recall that in the conformal class of $g$ there exists a \underbar{Gauduchon} metric $\tilde{g}$, i.e. a metric satisfying\ \ $\dbar\partial(\omega_{\tilde{g}}^{n-1}) = 0\ ;$ $\tilde{g}$ is unique up to a constant positive factor if\ \ \mb{n \geq 2}\ \ ([G] p. 502, [LT] Theorem 1.2.4).

There is a natural way to define a map
\[
\deg_g : \mathcal{H}(E) \map \mbb{R}\ ,
\]
called $g$-\underbar{degree}, with the following properties (see [LT] sections 1.3 and 1.4):\\
-\ If $g$ is a Gauduchon metric, and\ \ $\dbar \in \mathcal{H}(E)$\ \ is a holomorphic structure, then $\deg_g(\dbar)$ is given as follows: Choose any 
Hermitian metric $h$ in $E$, and let $d$ be the \underbar{Chern} 
\underbar{connection} in $(E,\partial)$ induced by $h$, i.e. the unique 
$h$-unitary connection in $E$ with $(0,1)$-part $\dbar$. Then 
\[ 
\deg_g(\dbar) := {i\over{2\pi}}\int\limits_X \tr(d^2)\wedge\omega_g^{n-1} =
{i\over{2n\pi}}\int\limits_X \tr\Lambda_g d^2\cdot\omega_g^n = 
{i\over{2n\pi}}\int\limits_X \tr\Lambda_g [\dbar,\partial]\cdot\omega_g^n\ .
\]
-\ If $g$ is arbitrary, then there is a unique Gauduchon metric $\tilde{g}$ in the conformal class of $g$ such that\ \ \mb{\deg_g = \deg_{\tilde{g}}\ .}

The $g$-\underbar{slope} of 
$\dbar$ is
\[
  \mu_g(\dbar) := {{\deg_g(\dbar)}\over r}\ ,
\]
where $r$ is the rank of $E$.

If\ \ $D = D^\prime + D\dpr$\ \ is a flat connection, then it holds\ \
$(D\dpr)^2 = 0\ ,$\ so $D\dpr$ is a holomorphic structure. We define the
$g$-degree and $g$-slope of $D$ as  
\[
\deg_g(D) := \deg_g(D\dpr)\ \ ,\ \ \mu_g(D) := \mu_g(D\dpr)\ .
\]
Similarly, for an integrable Higgs operator\ \ $d\dpr = \dbar + \theta$\ \
it holds\ \ $\dbar^2 = 0\ ,$\ and we define  
\[ 
\deg_g(d\dpr) := \deg_g(\dbar)\ \ ,\ \ \mu_g(d\dpr):=\mu_g(\dbar)\ .
\]
Observe that in all three cases the $g$-degrees (resp. slopes) of isomorphic operators are the same.

\begin{rem}\label{rem3}
Suppose that $g$ is a \underbar{K\"ahler} metric, i.e.\ \
$d(\omega_g) = 0\ .$ Then the $g$-degree is a topological
invariant of the bundle $E$, completely determined by the first real Chern
class\ \ $c_1(E)_\mbb{R} \in H^2(X,\mbb{R}).$\ In particular, since all real
Chern classes of a flat bundle vanish, it holds\ \ $\deg_g(D) =
0$\ \ for every flat connection $D$ in $E$. On the other hand, if e.g. 
$X$ is a surface admitting no K\"ahler metric and $g$ is Gauduchon, then
\underbar{every} real number is the $g$-degree of a flat line bundle on $X$ ([LT]
Proposition 1.3.13).
\end{rem}

\begin{lem}\label{lem2} If $g$ is a Gauduchon metric, then
for any metric $h$ in $E$ it holds:\\
i) If $D$ is a flat connection, then
\[ 
\deg_g(D) = -{i\over{n\pi}}\int\limits_X \tr\Lambda_g G_h\cdot\omega_g^n\ ,  
\]
where $G_h$ is the pseudocurvature of $d\dpr$ with respect to $h$.\\
ii) If $d\dpr$ is an integrable Higgs operator, then
\[ 
\deg_g(d\dpr) = {i\over{2n\pi}}\int\limits_X \tr\Lambda_g
F_h\cdot\omega_g^n\ ,  
\] 
where $F_h$ is the curvature of $d\dpr$ with respect to $h$.
\end{lem}

{\bf Proof:} i) Observe that\ \ $\Lambda_g G_h = \Lambda_g\dbar_h(\theta_h)\
.$\ \ The Chern connection in $(E,D\dpr)$ induced by $h$ 
is\ \ $D\dpr + \partial_h - \theta_h = D - 2\theta_h\ ,$\ and it holds 
\[
\tr\Lambda_g(D - 2\theta_h)^2 = -2\tr\Lambda_g((\dbar + \theta^*)(\theta)
= -2\tr\Lambda_g(G_h + [\theta,\theta^*]) = -2\tr\Lambda_g(G_h)\ ,
\]
so the claim follows by integration. \\ 
ii) Lemma~\ref{lem1} implies\ \ $\tr\Lambda_gF_h = \tr\Lambda_g d_h^2\ ;$\
again the claim follows by integration.  \qed

\section{Einstein metrics and stability for flat bundles.}

\newcommand{\vol}{\mathrm{vol}}

We fix a Hermitian metric $g$ in $X$; the associated volume form is\
\ $\vol_g := {1\over{n!}}\omega_g^n\ ,$\ and the $g$-volume of $X$ is\ \
$\Vol_g(X) := \int\limits_X \vol_g\ .$
We further fix a Hermitian metric $h$ in $E$, and denote by $\vert\ .\ \vert$
the pointwise norm on forms with values in $E$ (and associated bundles) defined
by $h$ and $g$.

Let\ \ $D \in \mathcal{A}_f(E)$\ \ be a flat connection in $E$, and write\ \ 
$D = d + \Theta = \partial + \dbar + \theta + \theta^*$\ \ as in section 1.
Let\ \ $d\dpr_h = I_h(D) = \dbar + \theta \in \mathcal{A}\dpr(E)$\ \ be the Higgs operator 
associated to  $D$, and\ \ $G_h = (d\dpr_h)^2$\ \ its pseudocurvature. From\
\ $\Lambda_g G_h = \Lambda_g\dbar_h(\theta_h)$\ \ and Lemma~\ref{lem1} we deduce 
\[ 
(i\Lambda_gG_h)^* = 
-i\Lambda_g((\dbar(\theta))^*) = -i\Lambda_g\partial(\theta^*) = 
i\Lambda_g\dbar(\theta) = i\Lambda_gG_h\ , 
\] 
so $i\Lambda_gG_h$ is selfadjoint with respect to $h$.

\begin{rem}\label{rem4}
It also holds\ \ 
$i\Lambda_gG_h = {i\over 2}\Lambda_g(\dbar(\Theta) - \partial(\Theta))\ 
,$\ which in the case of a \underbar{K\"ahler} metric $g$ equals $\half 
d^*(\Theta)$, where $d^*$ is the $L^2$-adjoint of\ \ $d = \partial + 
\dbar\ .$ 
\end{rem}

\begin{defn}
$h$ is called a $g$-\underbar{Einstein}
metric in  $(E,D)$ if\ \ $i\Lambda_gG_h = c\cdot\id_E$\ \ with a real
constant $c$, which is called the \underbar{Einstein} \underbar{constant}.
\end{defn}

\begin{lem}\label{exlem1}
Let $h$ be a $g$-Einstein metric in $(E,D)$, and\ \ \mb{\tilde{g} = \varphi\cdot g}\ \ conformally equivalent to $g$. Then there exists a $\tilde{g}$-Einstein metric $\tilde{h}$ in $(E,D)$ which is conformally equivalent to $h$.
\end{lem}

{\bf Proof:}\ \mb{\tilde{g} = \varphi\cdot g}\ \ implies\ \ \mb{\Lambda_{\tilde{g}} = \frac{1}{\varphi}\cdot\Lambda_g\ .}\ From Remark ~\ref{rem1} iii) it follows that for\ \ \mb{f \in \mathcal{C}^\infty(X,\mbb{R})}\ \ it holds\ \ \mb{G_{e^f\cdot h} = G_h - \frac{1}{4}\dbar\partial(f)\cdot\id_E\ .}\ Hence the condition\ \ \mb{i\Lambda_gG_h = c\cdot\id_E}\ \ implies\ \ \mb{i\Lambda_{\tilde{g}}G_{e^f\cdot h} = (\frac{c}{\varphi} - \frac{1}{4}P(f))\cdot\id_E\ ,}\ where\ \ \mb{P := i\Lambda_{\tilde{g}}\dbar\partial\ .}\ Since\ \ \mb{\mathcal{C}^\infty(X,\mbb{R}) = \im P\oplus\mbb{R}}\ ([LT] Corollary 2.9), there exists an $f$ such that $\frac{c}{\varphi} - \frac{1}{4}P(f)$ is constant. \qed

\begin{lem}\label{lem3}
If\ \ $i\Lambda_gG_h = c\cdot\id_E$\ \ with\ \ $c \in 
\mbb{R}\ ,$\ then it holds:\\
i)\ \mb{c = -{{\pi}\over{(n-1)!\cdot\Vol_g(X)}}\cdot\mu_g(D)}\ \ if $g$ is Gauduchon.\\
ii)\ $\deg_g(D) = 0$\ \ if and only if\ \ $c = 0\ .$
\end{lem}

{\bf Proof:} i) is an immediate consequence of Lemma~\ref{lem2}.\\
ii) If $g$ is Gauduchon, then this follows from i). If $g$ is arbitrary, then let\ \ \mb{\tilde{g} = \varphi\cdot g}\ \ be the Gauduchon metric in its conformal class such that\ \ \mb{\deg_g = \deg_{\tilde{g}}\ .} Now we have
\[
i\Lambda_gG_h = 0\ \ \Longleftrightarrow\ \ i\Lambda_{\tilde{g}}G_h = 0\ \ \Longleftrightarrow\ \ \deg_{\tilde{g}}(D) = 0\ \ \Longleftrightarrow\ \ \deg_g(D) = 0\ .
\]
\qed

\begin{rem}\label{rem5}
i) If two flat connections $D_1$,$D_2$ are isomorphic via the automorphism $f$ 
of $E$, i.e. if\ \ $D_2\circ f - f\circ D_1 = 0\ ,$\ and if $h$ is a
$g$-Einstein metric in $(E,D_1)$, then $f_*h$ is a $g$-Einstein metric
in $(E,D_2)$ with the same Einstein constant.\\
ii) By Remark~\ref{rem2}, a necessary condition for\ \ $d\dpr_h = I_h(D)$\ \ to
be an integrable  Higgs operator is that $h$ is a $g$-Einstein metric for $D$
with Einstein constant\ \ $c = 0\ ,$\  so in particular\ \ $\deg_g(D) = 0\ .$\ 
On the other hand it holds\ \ $d^2 = -\Theta\wedge\Theta$\ \ 
(Lemma~\ref{lem1}), and, if $d\dpr_h$ is  integrable,\ 
$\theta\wedge\theta = 0$\ \ implying\ \ $\theta^*\wedge\theta^* =  0\ 
.$\ This gives\ \ $\tr(d^2) = -\tr[\theta,\theta^*] = 0\ ,$\ which  
implies\ \ $\deg_g(d\dpr_h) = 0\ .$\\
iii) For complex vector bundles on compact \underbar{Riemannian} manifolds
$(X,g)$, Corlette defines a $g$-\underbar{harmonic} metric for a flat
connection by the condition\ \ $d^*(\Theta) =  0$\ \ ([C]). If $X$ is
complex and $g$ is a \underbar{K\"ahler} metric, then the $g$-degree of any
flat connection vanishes, so in this context $g$-harmonic is the same as
$g$-Einstein (see Remarks~\ref{rem3}\ and \ref{rem4}), but in general the 
two notions are different.
\end{rem}

Now we prove a useful Vanishing Theorem.

\begin{prop}\label{prop1}
Let $D$ be a flat connection in $E$, and $h$ a
$g$-Einstein  metric in $(E,D)$ with Einstein constant $c$. \\
If\ \ $c > 0\ ,$\ then the only section\ \ $s \in A^0(E)$\ \ with\ \ 
$D(s) = 0$\ \ is\ \ $s = 0\ .$ \\ 
If\ \ $c = 0\ ,$\ then for every section\ \ 
$s \in A^0(E)$\ \ with\ \ $D(s) = 0$\ \ it holds\ \ $\dbar(s) =
\theta(s) = 0$\ \ and\ \ $\partial(s) = \theta^*(s) = 0\ ,$\ so in
particular\ \ $d\dpr_h(s) = 0\  .$
\end{prop}

{\bf Proof:}\ $D(s) = 0$\ \ is equivalent to
\begin{equation}
\partial(s) = -\theta(s)\ \ ,\ \ \dbar(s) = -\theta^*(s)\ ; 
\end{equation}
this implies
\begin{equation}
\dbar\partial h(s,s) = -h(\dbar\circ\theta(s),s) - h(\theta(s),\theta(s)) + 
h(\dbar(s),\dbar(s)) - h(s,\partial\circ\theta^*(s))\ . 
\end{equation}
The assumption that $h$ is $g$-Einstein means\ \ $i\Lambda_g\dbar(\theta) 
= i\Lambda_gG_h = c\cdot\id_E\ ,$\ which is equivalent to
\ \ $i\Lambda_g\partial(\theta^*) = -c\cdot\id_E$\ \ since\ \ 
$(i\Lambda_g\dbar(\theta))^* = -i\Lambda_g(\dbar(\theta)^*) = 
-i\Lambda_g\partial(\theta^*)\ ;$\ these relations can be rewritten as
\begin{equation}
  i\Lambda_g\dbar\circ\theta = -i\Lambda_g\theta\circ\dbar + c\cdot\id_E\ \ ,\ 
\ i\Lambda_g\partial\circ\theta^* = -i\Lambda_g\theta^*\circ\partial - 
c\cdot\id_E\ . 
\end{equation}
Using $(1)$ and $(3)$ we get
\begin{eqnarray*}
  i\Lambda_gh(\dbar\circ\theta(s),s) &=&
-i\Lambda_gh(\theta\circ\dbar(s),s) + c\cdot\vert s\vert^2 = 
i\Lambda_gh(\dbar(s),\theta^*(s)) + c\cdot\vert s\vert^2 \\
&=& -i\Lambda_gh(\dbar(s),\dbar(s)) + c\cdot\vert s\vert^2 = 
\vert\dbar(s)\vert^2 + c\cdot\vert s\vert^2\ ,
\end{eqnarray*}
and similarly
\[
  i\Lambda_gh(s,\partial\circ\theta^*(s)) = \vert\theta(s)\vert^2 + 
c\cdot\vert s\vert^2\ ,
\]
so $(2)$ implies
\[
 i\Lambda_g\dbar\partial h(s,s) = -2\left(\vert\dbar(s)\vert^2 + 
\vert\theta(s)\vert^2 + c\cdot\vert s\vert^2\right)\ .
\]
Since the image of the operator $i\Lambda_g\dbar\partial$ on real functions
contains no non-zero functions of  constant sign ([LT] Lemma 7.2.7), this gives\
\ $s = 0$\ \ in the case\ \ $c > 0\ ,$\ \ and if\ \ $c = 0$\ \ we get\
\ $\dbar(s) = \theta(s) = 0\ ,$\ implying\ \ $\partial(s) = \theta^*(s) =
0$\ \ because of (1). \qed

The following corollary will be used later in the
context of moduli spaces.

\begin{cor}\label{cor1}
For\ \ $i = 1,2$\ \ let\ \ $D_i \in
\mathcal{A}_f(E)$\ \ be a flat connection, $h_i$ a $g$-Einstein metric
in $(E,D_i)$, and\ \ $d\dpr_i := I_{h_i}(D_i) \in \mathcal{A}\dpr(E)$\ \
the associated Higgs operator. If  $D_1$ and $D_2$ are
isomorphic via the automorphism $f$ of $E$, then $d\dpr_1$ and $d\dpr_2$ are
isomorphic via $f$, too.
\end{cor}

{\bf Proof:} Let $h$ be the metric in\ \ $\End E = E^*\otimes E$\ \ induced by
the dual metric of $h_1$ in $E^*$ and $h_2$ in $E$, and $D$ the connection in
$\End E$ defined by\ \ $D(f) = D_2\circ f - f\circ D_1$\ \ for all\ \ $f
\in A^0(\End E)\ .$\ Then $D$ is flat of $g$-degree 0 since $D_1$ and $D_2$
are flat of equal degree, and $h$ is a $g$-Einstein metric in $(\End E,D)$ with
Einstein constant\ \ $c = 0$\ (compare Remark~\ref{rem2}). Furthermore, the
Higgs operator $d\dpr$ in $\End E$ defined by\ \ $d\dpr(f) = d\dpr_2\circ f
- f\circ d\dpr_1$\ \ equals $I_h(D)$. Hence Proposition~\ref{prop1} implies that an
automorphism $f$ of $E$ with\ \ $D(f) = 0$\ \ also satisfies\ \
$d\dpr(f) = 0\ .$ \qed

If\ \ $F \subset E$\ \ is a $D$-invariant subbundle of $E$, then 
it is obvious that flatness of $D$ implies flatness of $D\vert_F$, and hence 
the following definition makes sense.

\begin{defn}
A flat connection $D$ in $E$ is called 
$g$-(\underbar{semi})\underbar{stable} iff for every proper 
$D$-invariant subbundle\ \ $0 \neq F \subset E$\ \ it holds\ \ 
$\mu_g(D\vert_F) > \mu_g(D)$\ ($\mu_g(D\vert_F) \geq \mu_g(D))$. 
$D$ is called $g$-\underbar{polystable} iff\ \ $E = E_1\oplus 
E_2\oplus\ldots\oplus E_k$\ \ is a direct sum of $D$-invariant and 
$g$-stable subbundles $E_i$ with\ \ $\mu_g(D\vert_{E_i}) = \mu_g(D)$\ \ for\
\ $i = 1,2,\ldots,k\ .$
\end{defn}

\begin{rem}\label{rem6}
i) Let $D$ be a flat connection in $E$, and\ \ $0 \neq F \subset E$\ \ 
a proper $D$-invariant subbundle. Then $g$-stability of $D$ implies\ \ \mb{\mu_g(D\vert_F) > \mu_g(D)}\ \ and hence the $g$-\underbar{instability} of the holomorphic 
structure $D\dpr$ in $E$ (in the sense of e.g. [LT]) since $F$ is a $D\dpr$-holomorphic subbundle of $E$.\\
ii) Suppose that $g$ is a K\"ahler metric; then\ \ 
$\deg_g(D) = 0$\ \ for every flat connection $D$ (Remark~\ref{rem3}). 
Hence a flat connection $D$ in $E$ is\\
\ -\ always $g$-semistable,\\
\ -\ $g$-stable if and only if $E$ has no proper non-trivial $D$-invariant 
subbundle,\\
\ -\ $g$-polystable if $E$ is a direct sum of $D$-invariant 
$g$-stable subbundles.\\
This means that $g$-(poly-)stability on a K\"ahler manifold coincides with 
(poly-)stability in the sense of Corlette [C].\\
iii) It is obvious that stability in the sense of Corlette always implies 
$g$-stability, but at the end of this section we will give an example of a 
$g$-stable bundle which is not stable in the sense of Corlette.
\end{rem}

\begin{defn}
A flat connection $D$ in $E$ is \underbar{simple} if
the only $D$-parallel endomorphisms $f$, i.e. those with
\ \ $D_\End(f) = D\circ f - f\circ D = 0\ ,$\ 
are the homotheties\ \ $f =
a\cdot\id_E\ ,\ a \in \mbb{C}\ .$
\end{defn}

Let $D$ be a flat connection in $E$,\ \ $0 \neq F \subset E$\ \ a
$D$-invariant subbundle, and\ \ $Q := \qmod{E}{F}$\ \ the quotient with
natural projection\ \ $\pi : E \map Q\ .$ Then $D$ induces a flat
connection $D_Q$ in $Q$ such that\ \ $D_Q\circ\pi = \pi\circ D\ .$\ In
particular, $F$ is a holomorphic subbundle of $(E,D\dpr)$, and $D_Q\dpr$
is the induced holomorphic structure in $Q$. Since the $g$-degree of a
flat connection $D$ by definition equals the $g$-degree of the
associated holomorphic structure $D\dpr$, it follow\ \ $\deg_g(D) =
\deg_g(D_1) + \deg_g(D_Q)\ .$\ Hence as in the case of holomorphic bundles
one verifies (compare [K] Chapter V)

\begin{prop}\label{prop2}
i) A flat connection $D$ in $E$ is $g$-(semi)stable if and
only if for every $D$-invariant proper subbundle\ \ $0 \neq F \subset E$\
\ with quotient\ \ $Q = \qmod{E}{F}$\ \ it holds\ \ $\mu_g(D_Q) <
\mu_g(D)$\ \ (resp. $\mu_g(D_Q) \leq \mu_g(D)$\ .)\\
ii) Let $(E_1,D_1)$ and $(E_2,D_2)$ be $g$-stable flat bundles over $X$ with\ \
$\mu_g(D_1) = \mu_g(D_2)\ .$\ If\ \ $f \in A^0(\Hom(E_1,E_2))$\ \
satisfies\ \ $D_2\circ f = f\circ D_1\ ,$\ then either\ \ $f = 0$\ \ or
$f$ is an isomorphism.\\
iii) A $g$-stable flat connection $D$ in $E$ is simple.
\end{prop}

Next we prove the first half of the main result of this section.

\begin{prop}\label{prop3}
Let $D$ be a flat connection in $E$, and $h$ a
$g$-Einstein  metric in $(E,D)$ with Einstein constant $c$; then $D$ is $g$-semistable. If $D$ is not
$g$-stable, then  $D$ is $g$-polystable; more precisely, \ $E = E_1\oplus 
E_2\oplus\ldots\oplus E_k$\ \ is a $h$-orthogonal direct sum of $D$-invariant 
$g$-stable subbundles such that\ \ $\mu_g(D\vert_{E_i}) = \mu_g(D)$\ \ for\ 
\ $i = 1,2,\ldots,k\ .$\ Furthermore, $h\vert_{E_i}$ is a $g$-Einstein
metric in $(E_i,D\vert_{E_i})$ with Einstein constant $c$ for all $i$, and the direct sum is invariant
with  respect to the Higgs operator\ \ $d\dpr_h = I_h(D)\ .$
\end{prop}

{\bf Proof:} First we consider the case when $g$ is a Gauduchon metric. Let\ \ $0 \neq F \subset E$\ \ be a $D$-invariant proper 
subbundle of rank $s$; then\ \ $E = F\oplus F^{\perp}\ ,$\ where $F^{\perp}$ is
the  $h$-orthogonal complement of $F$. With respect to this decomposition, we
write operators as $2\times 2$ matrices, so $D$ has the form 
\[ D = 
\left(
\begin{array}{cc}
D_1&A\\ 0&D_2
\end{array} 
\right)\ ,
\]  
where\ \ $D_1 = D\vert_F$\ \ and $D_2$ is a flat connection in $F^{\perp}$. 
We use notations as in section 2; it is easy to see that the 
operator $\delta$ associated to $D$ by $h$ has the form 
\[ \delta = \pmatrix{\delta_1&0\cr A^*&\delta_2}\ , \] 
where the $\delta_i$ are the operators associated to the $D_i$ by $h$. 
Similarly it holds
\[
  \dbar = \half(D\dpr + \delta\dpr) = \half
\left(
\begin{array}{cc}
D\dpr_1 + 
\delta\dpr_1&A\dpr\\ {A^\prime}^*&D\dpr_2 + \delta\dpr_2
\end{array} 
\right) = 
\left(
\begin{array}{cc}
\dbar_1&\half A\dpr\\ \half A{^\prime}^*&\dbar_2
\end{array} 
\right)\ , 
\]
and
\[
\theta = \half(D^\prime - \delta^\prime) = 
\left(
\begin{array}{cc}
D^\prime_1 - 
\delta^\prime_1&A^\prime\\ -{A\dpr}^*&D^\prime_2 - \delta^\prime_2
\end{array} 
\right) 
= 
\left(
\begin{array}{cc}
\theta_1&\half A^\prime\cr -\half{A\dpr}^*&\theta_2
\end{array} 
\right)\ , 
\]
where $A^\prime$ resp. $A\dpr$ is the part of $A$ of type $(1,0)$ resp. $(0,1)$. 
This implies
\begin{eqnarray*}
\dbar(\theta) &=& [\dbar,\theta]\\ &=& 
\left(
\begin{array}{c}
\dbar_1(\theta_1) + {1\over 4}(A^\prime\wedge{A^\prime}^*
- A\dpr\wedge{A\dpr}^*)\phantom{MMMMM}*\phantom{MMMMM}\\
\phantom{MMMMM}*\phantom{MMMMM}\dbar_2(\theta_2) + {1\over 4}({A^\prime}^*\wedge A^\prime -
{A\dpr}^*\wedge A\dpr)
\end{array} 
\right)\ , 
\end{eqnarray*}
hence
\begin{eqnarray}
c\cdot\id_E &=& i\Lambda_g G_h \nonumber\\
&=& 
\left(
\begin{array}{c}
i\Lambda_g G_{1,h} + {i\over 4}\Lambda_g(A^\prime\wedge{A^\prime}^*
- A\dpr\wedge{A\dpr}^*)\phantom{MMM}*\phantom{MMM}\\
\phantom{MMM}*\phantom{MMM}i\Lambda_g G_{2,h} + {i\over
4}\Lambda_g({A^\prime}^*\wedge A^\prime - {A\dpr}^*\wedge A\dpr)\end{array} 
\right)\ , 
\end{eqnarray}
and thus
\[
sc = \tr(i\Lambda_g G_{1,h} + {i\over 4}\Lambda_g(A^\prime\wedge{A^\prime}^*
- A\dpr\wedge{A\dpr}^*)) = i\tr\Lambda_g G_{1,h} + {1\over 4}\vert A\vert^2 \
. \]
Using Lemma~\ref{lem2} and Lemma~\ref{lem3} we conclude
\begin{equation}
\mu_g(D_1) = -{i\over{sn\pi}}\int\limits_X\tr\Lambda_g G_{1,h}\cdot\omega_g^n
\geq  -{{c(n-1)!}\over\pi}\Vol_g(X) = \mu_g(D)\ ;
\end{equation}
this prove that $D$ is $g$-semistable.\\
If $D$ is not $g$-stable, then there exists a subbundle $F$ as above such
that equality holds in $(5)$, which implies\ \ $A = 0\ .$\ This means
not only that $F^{\perp}$ is $D$-invariant, too, with\ \ $D\vert_{F^{\perp}} =
D_2$, but also that 
\[
i\Lambda_g G_{1,h} = c\cdot\id_F\ \ ,\ \ i\Lambda_g G_{2,h} =
c\cdot\id_{F^{\perp}} 
\]
by $(4)$. Hence the restriction of $h$ to $F$ resp. $F^{\perp}$ is
$g$-Einstein for $D_1$ resp. $D_2$, and it holds\ \ $\mu_g(D_1) = \mu_g(D) =
\mu_g(D_2)$\ \ by Lemma~\ref{lem3}. Furthermore, the $D$-invariance of $F$ means
that the inclusion\ \ $i : F \hookrightarrow E$\ \ is parallel with respect to
the flat connection in $\Hom(F,E)$ induced by $D_1$ and $D$. Using
Remark~\ref{rem2} and Proposition~\ref{prop1} as in the proof of 
Corollary~\ref{cor1}, we
conclude that $i$ is also parallel with respect to the associated Higgs
operator, i.e. that $F$ is $d\dpr_h$-invariant; the same argument works
for $F^{\perp}$. If $D_1$ and $D_2$ are stable, then we are done; otherwise
the proof is finished by induction on the rank. 

Now let $g$ be arbitrary, let $\tilde{g}$ be the Gauduchon metric in its conformal class with\ \ \mb{\deg_g = \deg_{\tilde{g}}\ ,}\ and let $\tilde{h}$ be a $\tilde{g}$-Einstein metric in the conformal class of $h$, which exists by Lemma~\ref{exlem1}; then the theorem holds for $\tilde{g}$ and $\tilde{h}$. Since $g$ and $\tilde{g}$ define the same degree and slope, and hence stability, it follows that $D$ is $\tilde{g}$-semistable. If $D$ is not $g$-stable, then there exists a $D$-invariant proper subbundle $F$ as above with\ \ \mb{\mu_{\tilde{g}}(D_1) = \mu_g(D_1) = \mu_g(D) = \mu_{\tilde{g}}(D)\ .}\ Note that the $h$-orthogonal complement $F^{\perp}$ of $F$ is also the $\tilde{h}$-orthogonal complement, since $h$ and $\tilde{h}$ are conformally equivalent. Hence, using $\tilde{g}$ and $\tilde{h}$ we conclude as above that\ \ \mb{D = 
\left(
\begin{array}{cc}
D_1&0\\ 0&D_2
\end{array} 
\right)}\ \ with respect the decomposition\ \ \mb{E = F\oplus F^{\perp}\ ;} now we can proceed as in the Gauduchon case.
\qed

Another consequence of Proposition~\ref{prop1} is

\begin{prop}\label{prop4} Let $D$ be a simple flat connection in $E$. If a
$g$-Einstein metric in $(E,D)$ exists, then it is unique up to a positive
scalar. \end{prop}

{\bf Proof:} Let $h_1$,$h_2$ be $g$-Einstein metrics in $(E,D)$, and\ \ $c
\in \mbb{R}$\ \ the Einstein constant. There are differentiable
automorphisms $f$ and $k$ of $E$, selfadjoint with respect to both $h_1$ and
$h_2$, such that\ \ $f = k^2$\ \ and\ \ $h_2(s,t) = h_1(f(s),t) =
h_1(k(s),k(t))$\ \ for all\ \ $s,t \in A^0(E)\ .$ Since $D$ is simple it
suffices to show\ \ $D(f) = 0\ .$

We define a new flat connection\ \ $\tilde{D} := k\circ D\circ k^{-1}\ .$
In what follows, operators $\delta$, $d$, $\Theta$ etc. with a subscript $i$
are associated to $D$ by the metric $h_i$, without a
subscript they are associated to $\tilde{D}$ by $h_1$. One verifies
\[
\delta_2 = f^{-1}\circ\delta_1\circ f\ \ ,\ \ \delta =
k^{-1}\circ\delta_1\circ k = k\circ\delta_2\circ k^{-1}\ , 
\]
implying
\[
d = \half(\tilde{D} + \delta) = k\circ d_2\circ k^{-1}\ \ ,
\ \ \Theta = \half(\tilde{D} - \delta) = k\circ\Theta_2\circ k^{-1}
\]
and hence
\[
i\Lambda_g G_{h_1} = i\Lambda_g\dbar(\theta) =
ik\circ\Lambda_g\dbar_2(\theta_2)\circ k^{-1} = ik\circ\Lambda_g G_{2,h_2}\circ
k^{-1} = c\cdot\id_E\ ,
\]
so $h_1$ is a $g$-Einstein metric in $(E,\tilde{D})$. It follows that $h_1$
induces a $g$-Einstein metric with Einstein constant 0 for the flat
connection\ \ $\tilde{D}_\End(.) = .\circ D - \tilde{D}\circ.$\ \ in $\End
E$. By definition it holds\ \ $\tilde{D}_\End(k) = 0\ ,$\ so Proposition~\ref{prop1}
implies\ \ $\tilde{d}_\End(k) = 0\ .$\ Since\ \ $\tilde{\delta}_\End =
2\tilde{d}_\End - \tilde{D}_\End\ ,$\ it follows
\[
0 = \tilde{\delta}_\End(k) = k\circ\delta_1 - \delta\circ k = k\circ\delta_1 -
k^{-1}\circ\delta_1\circ k^2 = k^{-1}\circ(f\circ\delta_1 - \delta_1\circ f)\
, 
\]
implying\ \ $\delta_{1,\End}(f) = 0\ ,$\ where $\delta_{1,\End}$ is the
operator on $\End E$ induced by $D$ and $h_1$. But this is equivalent to\ \
$\delta_{1,\End}^\prime(f) = 0$\ \ and\ \ $\delta\dpr_{1,\End}(f) = 0\
,$\ and taking adjoints with respect to $h_1$ we get
\[
0 = (\delta_{1,\End}^\prime(f))^* = D\dpr_\End(f)\ \ ,\ \ 0 =
(\delta_{1,\End}\dpr(f))^* = D^\prime_\End(f)\ , 
\]
i.e.\ \ $D_\End(f) = 0\ .$ \qed

Let $(E,D)$, $(\tilde{E},\tilde{D})$ be flat bundles with $g$-Einstein
metrics $h$, $\tilde{h}$. Let\ \ $E = \bigoplus\limits_{i=1}^k E_i$\ \ and\
\ $\tilde{E} = \bigoplus\limits_{i=1}^l \tilde{E}_i$\ \ be the orthogonal,
invariant splittings given by Proposition~\ref{prop3}. We write\ \ $D_i :=
D\vert_{E_i}\ ,$\ $\tilde{D}_i := \tilde{D}\vert_{\tilde{E}_i}\ ,$\ $h_i
:= h\vert_{E_i}\ ,$\ $\tilde{h}_i := \tilde{h}\vert_{\tilde{E}_i}\ .$\ Using
Propositions~\ref{prop2} and \ref{prop4} one verifies

\begin{cor}\label{cor2} Suppose that there exists an isomorphism\ \ $f \in
A^0(\Hom(E,\tilde{E}))$\ \ satisfying\ \ $f\circ D = \tilde{D}\circ f\ .$\ 
Then it holds\ \ $k = l\ ,$\ and, after renumbering of the summands if
necessary, there are isomorphisms\ \ $f_i \in
A^0(\Hom(E_i,\tilde{E}_i))$\ \ such that\ \ $f_i\circ D_i =
\tilde{D}_i\circ f$\ \ and\ \ $f_*(h_i) = \tilde{h}_i\ .$ \end{cor}

The following result is the converse of Proposition~\ref{prop3}.

\begin{prop}\label{prop5} Let $(E,D)$ a $g$-stable flat bundle over $X$.
Then there exists a $g$-Einstein metric for $(E,D)$. \end{prop}

{\bf Sketch of proof:} The proof is very similar to the one for the existence
of a $g$-Hermitian Einstein metric in a $g$-stable holomorphic vector bundle as
given in Chapter 3 of [LT]. Therefore we will be brief, leaving it to the
reader to fill in the necessary details.

First observe that by Lemma~\ref{exlem1} we may assume that $g$ is a Gauduchon metric.

For any metric $h$ in $E$ it holds
\[
G_h = \dbar_h(\theta_h) = 
{1\over 4}[D\dpr+\delta_h\dpr,D^\prime-\delta_h^\prime] 
= -{1\over 4}[D\dpr,\delta_h^\prime] + {1\over 4}[D^\prime,\delta_h\dpr]
\]
since\ \ $D^2 = \delta_h^2 = 0\ .$\ Observe that $[D\dpr,\delta_h^\prime]$
resp. $[D^\prime,\delta_h\dpr]$ is the curvature of the $h$-unitary
connection $D\dpr+\delta_h^\prime$ resp. $D^\prime+\delta_h\dpr$.

Fix a metric $h_0$ in $E$, and let\ \
$\delta = \delta^\prime + \delta\dpr\ ,$\ $d = \partial + \dbar\ ,$\
$\Theta = \theta + \theta^*$\ \ be the operators associated to\ \ $D =
D^\prime + D\dpr$\ \ and $h_0$ as in section 2. Consider for an
$h_0$-selfadjoint positive definite endomorphism $f$ of $E$ and\ \ $\eps \in
[0,1]$\ \ the differential equation     
\begin{equation}
L_\eps(f) := K^0 - {i\over 4}\Lambda_gD\dpr(f^{-1}\circ\delta^\prime(f)) + 
{i\over 4}\Lambda_gD^\prime(f^{-1}\circ\delta\dpr(f)) -\eps\cdot\log(f) = 0\ ,
\end{equation} 
where\ \
$K^0 := i\Lambda_g\dbar(\theta) - c\cdot\id_E = -{i\over
4}\Lambda_g([D\dpr,\delta^\prime] - [D^\prime,\delta\dpr]) - c\cdot\id_E\ ,$\
and $c$ is the constant associated to a possible $g$-Einstein metric
for $(E,D)$. The metric\ \ $f\cdot h_0\ ,$\ defined by\ \ $f\cdot
h_0(s,t) := h_0(f(s),t)$\ \ for sections $s,t$ in $E$, is $g$-Einstein
if and only if\ \ $L_0(f) = 0\ .$\ 

The term
\ \ $T_1 := i\Lambda_gD\dpr(f^{-1}\circ\delta^\prime(f))$\ \ (associated
to the unitary connection\ \ $d_1 := \delta^\prime + D\dpr$\ ) in equation
$(6)$ is of precisely the same type as the term\ \
\hbox{$T_0 := i\Lambda_g\dbar(f^{-1}\circ\partial_0(f))$}\ \ (associated to the
unitary connection\ \ $d_0 = \partial_0 + \dbar$\ ) in equation $(**)$ on
page 62 in [LT], and the term\ \ $T_2 := -
i\Lambda_gD^\prime(f^{-1}\circ\delta\dpr(f))$\ \ (associated to the unitary
connection\ \ $d_2 := D^\prime + \delta\dpr$\ ) is {\it almost} of this
type; e.g. the trace of all three terms equals
$i\Lambda_g\dbar\partial(\tr(\log f)$, and the symbols of the 
differential operators ${d\over{df}}\hat{T}_i$, where\ \ $\hat{T}_i(f) :=
f\circ T_i(f)\ ,$\ are equal, too. Therefore most of the arguments in [LT]
can easily be adapted to show first that for a simple flat 
connection $D$ equation $(6)$ has solutions $f_\eps$ for all\ \
$\eps \in (0,1]\ ,$\ which satisfy\ \ $\det f_\eps \equiv 1\ ,$\ 
and which converge to a solution $f$ of\ \ $L_0(f) = 0$\ \ if the 
$L^2$-norms of the $f_\eps$ are uniformly bounded. (There are two
places where one has to argue in a slightly different way: In the proof of the
analogue of [LT] Lemma 3.3.1, one uses the Laplacian\ \ $\Delta_D =
D^*\circ D$\ \ instead of $\Delta_{\dbar}$, and in the proof of the analogue of
[LT] Proposition 3.3.5 the \underbar{sum}\ \ $\Delta_{d_1} + \Delta_{d_2}$ of
the two Laplacians associated to $d_1$ and $d_2$ instead of just one.)

Then, under the assumptions that\ \ $\rk E \geq 2$\ \ and that the
$L^2$-norms of the $f_\eps$ are unbounded, one shows that for suitable \ \
$\eps_i \map 0\ ,$\ $\rho(\eps_i) \map 0\ ,$\ the limit
\[
\pi := \id_E - \lim_{\sigma\map
0}\left(\lim_{i\map\infty}\rho(\eps_i)\cdot f_{\eps_i}\right)^\sigma 
\]
exists weakly in $L^2_1$, and satisfies in $L^1$\ \ $\pi = \pi^* = \pi^2$\
\ and 
\begin{equation}
(\id_E - \pi)\circ D(\pi) = 0\ .
\end{equation}
This implies\ \ $(\id_E - \pi)\circ D\dpr(\pi) = 0\ ,$\ so
$\pi$ defines a weakly holomorphic subbundle $\mathcal{F}$ of the holomorphic bundle
$(E,D\dpr)$ by a theorem of Uhlenbeck and Yau (see [UY], [LT] Theorem 3.4.3).
$\mathcal{F}$ is a coherent subsheaf of $(E,D\dpr)$, a holomorphic subbundle
outside an analytic subset\ \ $S \subset X$\ \ of codimension at least 2,
and $\pi$ is smooth on $X\setminus S$. Therefore $(7)$ implies that
$\mathcal{F}\vert_{X\setminus S}$ is in fact a $D$-invariant subbundle of
$E\vert_{X\setminus S}$, which extends to a $D$-invariant subbundle $F$ of $E$
by the Lemma below. Again using arguments as is [LT], one finally shows that
$F$ violates the stability condition for $(E,D)$. \qed

\begin{lem} Let $X$ be a differentiable manifold, $E$ a differentiable
vector bundle over $X$, and $D$ a flat connection in $E$. Let\ \ $S \subset
X$\ \ be a subset such that $X\setminus S$ is open and dense in $X$, and with
the following property: For every point\ \ $x \in S$\ \ and every open
neighborhood $U$ of $x$ in $X$ there exists an open neighborhood\ \ $x \in
U^\prime \subset U$\ \ such that $U^\prime\setminus S$ is path-connected.\\
Then every $D$-invariant subbundle $\mathcal{F}$ of $E\vert_{X\setminus S}$
extends to a $D$-invariant subbundle $F$ of $E$.\end{lem}

{\bf Proof:} For every\ \ \mb{x \in S}\ \ choose an open neighborhood\ \ \mb{x
\in U \subset X}\ \ such that $U\setminus S$ is path connected and\ \
\mb{(E\vert_U,D) \cong (U\times V,d)\ ,}\ where $V$ is a vector space and $d$
the trivial flat connection. Since $\mathcal{F}$ is $D$-invariant and $U\setminus S$ 
is path connected, it holds
\[
(\mathcal{F}\vert_{U\setminus S},D) \cong ((U\setminus S)\times W,d)\ ,
\]
where\ \ \mb{W \subset V}\ \ is a constant subspace. Define $F$ over $U$ by\ \
\mb{F\vert_U :\cong U\times W\ ;}\ then the topological condition on $S$
implies that this is well defined on $S$, and hence gives a $D$-invariant
extension $F$ of $\mathcal{F}$ over $X$. \qed

The following main result of this section is a direct consequence of
Propositions~\ref{prop3} and \ref{prop5}.

\begin{thm}\label{thm1} A flat connection $D$ in $E$ admits a $g$-Einstein
metric if and only if it is $g$-polystable. \end{thm}

As for stable vector bundles and Hermitian-Einstein metrics, the gauge
theoretic interpretation of our results is as follows. The group
\[
\mathcal{G}^\mbb{C} := A^0(GL(E))
\]
of differentiable automorphisms of $E$ acts on $\mathcal{A}(E)$ by\ \
\mb{D\cdot f = f^{-1}\circ D\circ f\ ,}\ so
\[
\qmod{\mathcal{A}(E)}{\mathcal{G}^\mbb{C}}
\]
is the moduli space of isomorphism classes of connections in $E$. Observe that
flatness, simplicity and $g$-stability are preserved under this action. Fix a
metric $h$ in $E$; then it holds:

\begin{cor} The following two statements for a flat
connection $D$ are equivalent:\\
i) $D$ is $g$-stable.\\
ii) $D$ is simple, and there is a connection $D_0$ in the $\mathcal{G}^\mbb{C}$-orbit
through $D$ such that $h$ is $g$-Einstein for $D_0$.
\end{cor}

The essential uniqueness of a $g$-Einstein metric (Proposition~\ref{prop4}) implies
that the connection $D_0$ in ii) is unique up to the action of the subgroup
\[
\mathcal{G} := A^0(U(E,h)) \subset \mathcal{G}^\mbb{C}
\]
of $h$-unitary automorphisms. This
means that the moduli space
$$
\mathcal{M}_f^{\rm st}(E) = \qmod{\sett{D \in \mathcal{A}_f(E)}{D\ \mbox{is}\ g-\mbox{
stable}}}{\mathcal{G}^\mbb{C}} 
$$
of isomorphism classes of $g$-stable flat connections in $E$ coincides with the
quotient
$$
\qmod{\sett{D \in \mathcal{A}_f(E)}{D\ \mbox{is
simple and}\ h\ \mbox{is}\ g-\mbox{Einstein for}\ D}}{\mathcal{G}}\ .
$$

\vspace{.5cm}

{\bf Example:} We now give the promised example of a flat bundle which is 
$g$-stable, but not stable in the sense of Corlette.

An Inoue surface of type $S^\pm_N$ is the quotient of $\mbb{H}\times\mbb{C}$ by 
an affine transformation group $G$ generated by 
\begin{eqnarray*}
g_0(w,z) &:=& (\alpha w,\pm z+t)\ ,\\
g_i(w,z) &:=& (w+a_i,z+b_iw+c_i)\ ,\ i = 1,2,\\
g_3(w,z) &:=& (w,z+c_3)\ ,
\end{eqnarray*}
with certain constants\ \ \mb{\alpha,a_i,b_i,c_3 \in \mbb{R}\ ,}
\ \mb{c_1,c_2 \in \mbb{C}}\ \ (see [P] p. 160). Since the second 
Betti number of $S^\pm_N$ vanishes, the degree map
\[
\deg_g : \mbox{Pic}(S^\pm_N) \map \mbb{R}
\]
associated to a Gauduchon metric $g$ is, up to a positive factor, 
independent of the chosen metric $g$. In particular, \underbar{all} Hermitian metrics $g$ 
define the same notion of $g$-stability ([LT] Remark 1.4.4 iii)). 

The trivial flat connection $d$ on $\mbb{H}\times\mbb{C}$ induces a 
flat connection $D$ in the tangent bundle\ \ \mb{E := T_{S^\pm_N}\ .}\ 
A $D$-invariant sub-line bundle of $E$ is in particular a holomorphic 
subbundle, so it defines a holomorphic foliation of $S^\pm_N$. 
According to [B] Th\'eor\`eme 2, there is precisely one such foliation, 
namely the one induced by the $G$-invariant vertical foliation 
(i.e. with leaves $\{w\}\times\mbb{C}$) of $\mbb{H}\times\mbb{C}$. 
The corresponding trivial line bundle $L_0$ on $\mbb{H}\times\mbb{C}$ 
is $d$-invariant, so it descends to a unique $D$-invariant subbundle 
$L$ of $E$; this shows that $E$ is not stable in the sense of Corlette. 
Observe that $L$ has factors of automorphy
\ \ \mb{\chi(g_i) = \pm1\ ,\ i = 0,1,2,3\ ,}\ so the standard flat 
metric in $L_0$ defines a metric $h$ in $L$ such that the associated 
Chern connection in $(L,D\dpr\vert_L)$ is flat; this implies
\ \ \mb{\mu_g(D\vert_L) = \deg_g(D\vert_L) = 0\ .}\ On the other 
hand, the $g$-degree, and hence the $g$-slope, of $E$ is negative 
by [P] Proposition 4.7; this implies the $g$-stability of $E$ since 
$L$ is the only $D$-invariant proper subbundle of $E$.

\section{Einstein metrics and stability for Higgs bundles.}

Again we fix Hermitian metrics $g$ in $X$ and $h$ in $E$.

Let\ \ \mb{d\dpr = \dbar + \theta \in \mathcal{A}\dpr_i(E)}\ \ be an integrable 
Higgs operator,
\[
D_h = I_h^{-1}(d\dpr) = d + \Theta = \partial + \dbar + \theta + \theta^* \in 
\mathcal{A}(E)
\]
the connection associated to $d\dpr$ as in section 2, and\ \ \mb{F_h =
D_h^2}\ \ its curvature. 

\begin{defn} $h$ is called a $g$-\underbar{Einstein} metric in 
$(E,d\dpr)$ if\ \ \mb{K_h := i\Lambda_gF_h = c\cdot\id_E}\ \ with a real 
constant $c$, the \underbar{Einstein} \underbar{constant}. \end{defn}

\begin{lem}\label{exlem2}
Let $h$ be a $g$-Einstein metric in $(E,d\dpr)$, and\ \ \mb{\tilde{g} = \varphi\cdot g}\ \ conformally equivalent to $g$. Then there exists a $\tilde{g}$-Einstein metric $\tilde{h}$ in $(E,d\dpr)$ which is conformally equivalent to $h$.
\end{lem}

{\bf Proof:} From Remark ~\ref{rem1} iii) it follows that for\ \ \mb{f \in \mathcal{C}^\infty(X,\mbb{R})}\ \ it holds\ \ \mb{F_{e^f\cdot h} = F_h + \dbar\partial(f)\cdot\id_E\ .}\ Using this, the proof is analogous to that of Lemma~\ref{exlem1}.\qed

Notice that since $d\dpr$ is integrable it holds (compare Lemma~\ref{lem1})
\[
K_h = i\Lambda_g(d^2 + [\theta,\theta^*]) = i\Lambda_g([\partial,\dbar] +
[\theta,\theta^*])
\]
where\ \ \mb{d = \partial + \dbar\ .}  An immediate consequence of 
Lemma~\ref{lem2} and Lemma~\ref{exlem2} is (compare the proof of Lemma~\ref{lem3})

\begin{lem}\label{lem4}
If\ \ $i\Lambda_gF_h = c\cdot\id_E$\ \ with\ \ $c \in 
\mbb{R}\ ,$\ then it holds:\\
i)\ \mb{c = {{2\pi}\over{(n-1)!\cdot\Vol_g(X)}}\cdot\mu_g(d\dpr)}\ \ if $g$ is Gauduchon.\\
ii)\ $\deg_g(d\dpr) = 0$\ \ if and only if\ \ $c = 0\ .$
\end{lem}

\begin{rem}\label{rem7} (compare Remark \ref{rem5})\\
i) If two integrable Higgs operators $d\dpr_1$,$d\dpr_2$ are isomorphic via
the automorphism $f$  of $E$, i.e. if\ \ \mb{d\dpr_2\circ f - f\circ d\dpr_1 =
0\ ,}\ and if $h$ is a $g$-Einstein  metric in $(E,d\dpr_1)$, then $f_*h$ is a
$g$-Einstein metric in $(E,d\dpr_2)$, and the  associated Einstein
constants are equal.\\
ii) By Remark~\ref{rem2}, a necessary 
condition for\ \ \mb{D_h = I_h(d\dpr)}\ \ to be a flat connection is $h$ to be 
Einstein with Einstein constant\ \ \mb{c = 0\ ,}\ so in particular\ \ 
\mb{\deg_g(d\dpr) = 0\ .}\ On the other hand, the Chern connection in 
$(E,D\dpr_h)$ is $\partial -\theta + \dbar + \theta^*$, so the $g$-degree of 
$D_h$ is obtained by integrating $\tr\Lambda_g[\dbar + \theta^*,\partial - 
\theta]$ which equals $\tr\Lambda_g[\dbar,\partial]$ since $d\dpr$ is 
integrable (Lemma~\ref{lem1} ii)). If $D_h$ is flat, we furthermore have\ \ \mb{d^2 = 
-\Theta\wedge\Theta}\ (Lemma~\ref{lem1} i)), implying\ \ \mb{\tr[\dbar,\partial] = 0} 
and hence\ \ \mb{\deg_g(D_h) = 0\ .}\end{rem}

In analogy with the case of Hermitian-Einstein metrics in holomorphic vector
bundles, the following vanishing theorem holds.

\begin{prop}\label{prop6} Let $h$ be a $g$-Einstein metric in $(E,d\dpr)$ with
Einstein constant $c$. \\
If\ \ \mb{c < 0\ ,}\ then the only section\ \ \mb{s \in A^0(E)}\ \ with\ \
\mb{d\dpr(s) = 0}\ \ is\ \ \mb{s = 0\ .}\\
If\ \ \mb{c = 0\ ,}\ then for every section\ \ \mb{s \in A^0(E)}\ \ with\ \
\mb{d\dpr(s) = 0}\ \ it holds\ \ \mb{D_h(s) = 0\ .} \end{prop}

{\bf Proof:} For\ \ \mb{s \in A^0(E)\ ,}\ \mb{d\dpr(s) = 0}\
\ is equivalent to\ \ \mb{\dbar(s) = 0 = \theta(s)\ .}\ This implies
\begin{equation}
c\cdot\vert s\vert^2 = c\cdot h(s,s) = h(K_h(s),s) =
i\Lambda_g\left(h(\dbar\partial(s),s) + h(\theta^*(s),\theta^*(s))\right)\ . 
\end{equation}
We have
\[
i\Lambda_g\dbar\partial h(s,s) =
i\Lambda_g\left(h(\dbar\partial(s),s) - h(\partial(s),\partial(s))\right)
\]
since\ \ \mb{\dbar(s) = 0\ ,}\ and using $(8)$ we get 
\[
i\Lambda_g\dbar\partial h(s,s) = c\cdot\vert s\vert^2 - \vert\partial(s)\vert^2
- \vert\theta^*(s)\vert^2 \ . 
\]
Now the claim follows as in the proof of Proposition~\ref{prop1}. \qed

The proof of the following corollary is analogous to that of Corollary~\ref{cor1}.

\begin{cor}\label{cor3} For\ \ \mb{i = 1,2}\ \ let\ \ \mb{d\dpr_i \in
\mathcal{A}\dpr_i(E)}\ \ be an integrable Higgs operators, $h_i$ a $g$-Einstein metric
in $(E,d\dpr_i)$, and\ \ \mb{D_i := I_{h_i}^{-1}(d\dpr_i) \in \mathcal{A}(E)}\ \
the associated connection. If $d\dpr_1$ and $d\dpr_2$ are
isomorphic via the automorphism $f$ of $E$, then $D_1$ and $D_2$ are
isomorphic via $f$, too.\end{cor}

Let\ \ \mb{d\dpr = \dbar + \theta}\ \ be an integrable Higgs operator in $E$. A
coherent subsheaf $\mathcal{F}$ of the holomorphic bundle $(E,\dbar)$ is called a
\underbar{Higgs}-\underbar{subsheaf} of $(E,d\dpr)$ iff it is
$d\dpr$-invariant. For the definition of the $g$-degree and $g$-slope of a
coherent sheaf see [LT].

\begin{defn} An integrable Higgs operator $d\dpr$ in $E$ is called
$g$-(\underbar{semi})stable iff for every coherent Higgs-subsheaf $\mathcal{F}$ of
$(E,d\dpr)$ with\ \ \mb{0 < \rk\mathcal{F} < \rk E}\ \ it holds\ \ \mb{\mu_g(\mathcal{F}) <
\mu_g(E)}\ (\ \mb{\mu_g(\mathcal{F}) \leq \mu_g(E)}\ ). $d\dpr$ is called
$g$-\underbar{polystable} iff $E$ is a direct sum\ \ \mb{E = E_1\oplus E_2\oplus\ldots\oplus E_k}\
\ of $d\dpr$-invariant and $g$-stable subbundles $E_i$ with\ \
\mb{\mu_g(d\dpr\vert_{E_i} = \mu_g(d\dpr)}\ \ for\ \ \mb{i = 1,2,\ldots,k\ .}
\end{defn}

\begin{defn} An integrable Higgs operator $d\dpr$ in $E$ is called
\underbar{simple} iff for every\ \ \mb{f \in A^0(\End E)}\ \ with\
\ \mb{d\dpr\circ f = f\circ d\dpr}\ \ it holds\ \ \mb{f =
a\cdot\id_E}\ \ with\ \ \mb{a \in \mbb{C}\ .}\end{defn}

As in the case of stable vector bundles or flat connections,
(semi)-stability can equivalently be defined using quotients of $E$; again it
follows  

\begin{lem}\label{lem5} 
i) A $g$-stable integrable Higgs operator in $E$ is simple.\\
ii) Let $d\dpr_1$, $d\dpr_2$ be $g$-stable integrable Higgs operators in
bundles $E_1$, $E_2$ on $X$ such that\ \ \mb{\mu_g(d\dpr_1) = \mu_g(d\dpr_2)\
.}\ If\ \ \mb{f \in A^0(\Hom(E_1,E_2))}\ \ satisfies\ \ \mb{d\dpr_2\circ f =
f\circ d\dpr_1\ ,}\ then either\ \ \mb{f = 0}\ \ or $f$ is an isomorphism.
\end{lem}

Furthermore, using arguments similar to those in the proof of 
Proposition~\ref{prop4}, we get the following consequence 
of Proposition~\ref{prop6}.

\begin{prop}\label{prop7} Let $d\dpr$ be a simple integrable Higgs operator in
$E$. If a $g$-Einstein metric in $(E,d\dpr)$ exists, then it is unique up to
a positive scalar. \end{prop}

The proof of the next result is a straightforward generalization of that in the
K\"ahler case [S2] (just as for the proof of the corresponding statement for
Hermite-Einstein metrics in vector bundles, see [LT]).

\begin{prop}\label{prop8} Let $d\dpr$ be an integrable Higgs operator in $E$, and
$h$ a $g$-Einstein metric in $(E,d\dpr)$ with Einstein constant $c$; then $d\dpr$ is $g$-semistable. If
$d\dpr$ is not $g$-stable, then  $d\dpr$ is $g$-polystable; more precisely, \
\mb{E = E_1\oplus  E_2\oplus\ldots\oplus E_k}\ \ is an $h$-orthogonal direct
sum of $d\dpr$-invariant and $g$-stable subbundles such that\ \
\mb{\mu_g(d\dpr\vert_{E_i}) = \mu_g(d\dpr)}\ \ for\  \ \mb{i = 1,2,\ldots,k\
.}\ Furthermore, $h\vert_{E_i}$ is a $g$-Einstein metric in
$(E_i,d\dpr\vert_{E_i})$ with Einstein constant $c$ for all $i$, and the direct sum is invariant with 
respect to the connection\ \ \mb{D_h = I_h^{-1}(d\dpr)\ .} \end{prop}

Let $d\dpr$, $\tilde{d}\dpr$ be integrable  Higgs operators in bundles
$E$, $\tilde{E}$ with $g$-Einstein metrics $h$, $\tilde{h}$. Let\ \ \mb{E
= \bigoplus\limits_{i=1}^k E_i}\ \ and\ \ \mb{\tilde{E} =
\bigoplus\limits_{i=1}^l \tilde{E}_i}\ \ be the orthogonal, invariant
splittings given by Proposition~\ref{prop8}. We write\ \ \mb{d\dpr_i := d\dpr\vert_{E_i}\
,}\ \mb{\tilde{d}\dpr_i := \tilde{d}\dpr\vert_{\tilde{E}_i}\ ,}\ \mb{h_i :=
h\vert_{E_i}\ ,}\ \mb{\tilde{h}_i := \tilde{h}\vert_{\tilde{E}_i}\ .}

As in the previous section (but now using Lemma~\ref{lem5} and Proposition~\ref{prop7}) we
deduce

\begin{cor}\label{cor4} Suppose that there exists an isomorphism\ \ \mb{f \in
A^0(\Hom(E,\tilde{E}))}\ \ satisfying\ \ \mb{f\circ d\dpr = \tilde{d}\dpr\circ
f\ .}\  Then it holds\ \ \mb{k = l\ ,}\ and, after renumbering of the summands
if necessary, there are isomorphisms\ \ \mb{f_i \in
A^0(\Hom(E_i,\tilde{E}_i))}\ \ such that\ \ \mb{f_i\circ d\dpr_i =
\tilde{d}\dpr_i\circ f}\ \ and\ \ \mb{f_*(h_i) = \tilde{h}_i\ .} \end{cor}

\begin{rem} We expect that the existence of a $g$-Einstein metric for a
$g$-stable Higgs operator $d\dpr$ can be proved by solving (again using the
continuity method as in [LT]) the differential equation
$$
K_h + i\Lambda_gd\dpr(f^{-1}\circ d^\prime(f)) = c\cdot\id_E
$$
for a positive definite and $h$-selfadjoint endomorphism $f$ of $E$, where $h$
is a suitable fixed metric in $E$. \end{rem}

\section{Surfaces.}

In this section we consider the special case\ \ \mb{n = 2\ ,}\ i.e. where $X$
is a compact complex surface; again we fix a Hermitian metric $g$ in $X$.
In this case, the real Chern numbers\ \ \mb{c_1^2(E),c_2(E) \in H^4(X,\mbb{R})
\cong \mbb{R}}\ \ can be calculated by integrating the corresponding Chern forms
of any connection in $E$, independently of the chosen metric $g$. In
particular, if $E$ admits a flat connection, then these Chern numbers vanish.  

\begin{prop}\label{prop9} Suppose that\ \ \mb{D \in \mathcal{A}_f(E)}\ \ is a flat
connection of $g$-degree 0, and that $h$ is a $g$-Einstein metric in $(E,D)$.
Then it holds\ \ \mb{G_h = 0\ .}\ In particular, the Higgs operator
$d\dpr_h$ associated to $D$ and $h$ is integrable with\ \ \mb{\deg_g(d\dpr_h) =
0\ ,}\ and $h$ is a $g$-Einstein metric for $(E,d\dpr_h)$. \end{prop}

{\bf Proof:} (see [S2]) For\ \ \mb{\epsilon > 0}\ \ we define a new connection\ \
\mb{B_\epsilon := d + {1\over\epsilon}\theta + \epsilon\theta^*\ ,}\ and\ \
\mb{F_\epsilon := B_\epsilon^2\ .}\ Observe that\ \ \mb{n = 2}\ \ implies\ \
\mb{F_\epsilon^2 = {1\over{\epsilon^2}}\nabla_\epsilon^4\ ,}\ where\ \
\mb{\nabla_\epsilon = d\dpr_h + \epsilon d^\prime\ .}\ The vanishing of the
Chern numbers of $E$ implies\ \ \mb{\int\limits_X \tr F_\epsilon^2 = 0
\ ,}\ and hence\ \ \mb{\int\limits_X \tr \nabla_\epsilon^4 = 0}\ \ for all\ \
\mb{\epsilon > 0\ .}\ Taking the limit\ \ \mb{\epsilon \rightarrow 0}\ \ it
follows
\begin{equation}
\int\limits_X \tr G_h^2 = 0\ .
\end{equation}
Write\ \ \mb{G_h = G_{1,1} + G_2\ ,}\ where $G_{1,1}$ is the component of the
2-form $G_h$ of type $(1,1)$. Then it holds
\begin{equation}
*_gG_{1,1} = -G_{1,1}\ \ ,\ \ *_gG_2 = G_2\ ; 
\end{equation}
the first equation is a consequence of\ \ \mb{\Lambda_gG_h = 0\ ,}\ which follows from the assumption and Lemma~\ref{lem3}.
On the other hand, it holds\ \ \mb{G_h = \dbar^2 + \dbar(\theta) +
\theta\wedge\theta\ ,}\ so Lemma~\ref{lem1} implies 
\begin{equation}
G_{1,1} = \dbar(\theta) = \partial(\theta^*)^* = -\dbar(\theta)^* =
-G_{1,1}^*\ ,
\end{equation}
and
\begin{equation}
G_2 = \dbar^2 + \theta\wedge\theta =
-\theta^*\wedge\theta^* - \theta\wedge\theta = (\theta\wedge\theta +
\theta^*\wedge\theta^*)^* = G_2^*\ .
\end{equation}
(11) and (12) combined with (10) give\ \ \mb{*_gG_h^* = G_h\ ,}\ so from (9) it
follows
$$
0 = \int\limits_X \tr G_h^2 = \int\limits_X \tr(G_h\wedge *_gG_h^*) =
\int\limits_X\vert G_h\vert^2\vol_g\ , 
$$
implying\ \ \mb{(d\dpr_h)^2 = G_h = 0\ .}\ Hence $d\dpr_h$ is integrable,
$\deg_g(d\dpr_h)$ vanishes (Remark~\ref{rem5}), and $h$ is $g$-Einstein for
$(E,d\dpr)h)$ because the curvature of $d\dpr_h$ with respect to $h$ equals\ \
\mb{D^2 = 0\ .}\qed

\begin{prop}\label{prop10} Suppose that\ \ \mb{c^2_1(E) = c_2(E) = 0\ ,}\ that
$d\dpr$ is an integrable Higgs operator of $g$-degree 0, and that $h$ is a
$g$-Einstein metric in $(E,d\dpr)$. Then it holds\ \ \mb{F_h = 0\ .}\ In
particular, the connection $D_h$ associated to $d\dpr$ and $h$ is flat with\
\ \mb{\deg_g(D_h) = 0\ ,}\ and $h$ is a $g$-Einstein metric for
$(E,D_h)$.  \end{prop}

{\bf Proof:} Define\ \ \mb{F_{1,1} := d^2 + [\theta,\theta^*]\ ,}\ \mb{F_2 :=
\partial(\theta) + \dbar(\theta^*)\ ;}\ then\ \ \mb{F_h = F_{1,1} + F_2\ .}\ 
Observe that $F_{1,1}$ is of type
(1,1) because $d$ is a unitary connection in the holomorphic
bundle $(E,\dbar)$. Since\ \ \mb{\deg_g(d\dpr) = 0\ ,}\ Lemma~\ref{lem4} implies\ \
\mb{0 = \Lambda_g F_h = \Lambda_g F_{1,1}\ ,}\ hence it holds\ \
\mb{*_gF_{1,1} = -F_{1,1}}\ \ and\ \ \mb{*_gF_2 = F_2\ .}\ On the other hand,
it is easy to see that\ \ \mb{F_{1,1}^* = -F_{1,1}}\ \ and\ \ \mb{F_2^* = F_2\
.}\ Combining these relations we get\ \ \mb{*_gF_h^* = F_h\ .}\ Since
$c_1^2(E)$ resp. $c_2(E)$ are obtained by integrating $-{{1}\over
{4\pi^2}}(\tr F_h)^2$ resp. $-{{1}\over {8\pi^2}}((\tr F_h)^2 - \tr(F_h^2))$,
we get
$$
0 = \int\limits_X \tr(F_h^2) = \int\limits_X \tr(F_h\wedge *_gF_h^*) = \Vert
F_h\Vert^2\ , 
$$
implying\ \ \mb{D_h^2 = F_h = 0\ .}\ Hence $D_h$ is flat, 
$\deg_g(D_h)$ vanishes (Remark~\ref{rem7}), and $h$ is $g$-Einstein for $(E,D_h)$
because the pseudocurvature of $D_h$ with respect to $h$ equals\ \
\mb{(d\dpr)^2 = 0\ .}\qed

\begin{rem}\label{rem9} The above proposition implies in particular the following:
Suppose that\ \ \mb{c_1^2(E) = c_2(E) = 0\ ;}\ if there exists an
integrable Higgs operator $d\dpr$ in $E$ with $g$-degree 0 admitting a
$g$-Einstein metric, then the real Chern \underbar{class}\ \ \mb{c_1(E)_\mbb{R}
\in H^2(X,\mbb{R})}\ \ vanishes, because there is a flat connection in $E$. 
\end{rem}

We define $\mathcal{A}_f(E)^0_g$ to be the space of\ \ \mb{D \in \mathcal{A}_f(E)}\ \ of
$g$-degree 0 such that there exists a $g$-Einstein metric in $(E,D)$, and
$\mathcal{A}\dpr_i(E)^0_g$ to be the space of\ \ \mb{d\dpr \in \mathcal{A}\dpr_i(E)}\ \ of
$g$-degree 0 such that there exists a $g$-Einstein metric in $(E,d\dpr)$. 
By Remark~\ref{rem5} and Remark~\ref{rem7}, the two moduli sets 
\[
\mathcal{M}_f(E)^0_g := \qmod{\mathcal{A}_f(E)^0_g}{\mbox{ isomorphy of connections}}
\]
and
\[
\mathcal{M}\dpr(E)^0_g := \qmod{\mathcal{A}\dpr_i(E)^0_g}{\mbox{isomorphy of Higgs operators}}
\]
are well defined. The main result of this section is

\begin{thm}\label{thm3} There is a natural bijection
$$
I : \mathcal{M}_f(E)^0_g \map \mathcal{M}\dpr(E)^0_g\ .
$$
\end{thm}

{\bf Proof:} First observe that we may assume that the real Chern classes
of $E$ vanish, since otherwise both spaces are empty (see Remark~\ref{rem9}).

Let $D$ be a flat connection in $E$ with $g$-degree 0, and
$h$ a $g$-Einstein metric in $(E,D)$. By Proposition~\ref{prop9}, the associated Higgs
operator\ \ \mb{d\dpr_h = I_h(D)}\ \ is integrable with $g$-degree 0, and $h$
is a $g$-Einstein metric in $(E,d\dpr_h)$. We will show that the map $I$
defined by\ \ \mb{I([D]) := [d\dpr_h]}\ \ is well defined and bijective. 

Suppose that\ \ \mb{D,\tilde{D} \in \mathcal{A}_f(E)^0_g}\ \ are isomorphic via
the automorphism $f$ of $E$; then $f_*h$ is $g$-Einstein in $(E,\tilde{D})$
(Remark~\ref{rem5}), the Higgs-operator $\tilde{d}\dpr$ associated to $\tilde{D}$
and $f_*h$ is isomorphic to $d\dpr$ via $f$ (Corollary~\ref{cor1}), and $f_*h$
is a $g$-Einstein metric in $(E,\tilde{d}\dpr)$ (Remark~\ref{rem7}). To prove that
$I$ is well defined it thus suffices to show that two different
$g$-Einstein metrics $h,\tilde{h}$ for a fixed\ \ \mb{D \in \mathcal{A}_f(E)^0_g}\
\ produce isomorphic Higgs operators $d\dpr_h,d\dpr_{\tilde{h}}$. For this
consider the $D$-invariant and $h$- resp. $\tilde{h}$-orthogonal splittings\ \
\mb{E = \bigoplus\limits_{i=1}^k E_i}\ \ resp.\ \ \mb{E =
\bigoplus\limits_{i=1}^l \tilde{E}_i}\ \ associated to $h$ resp. $\tilde{h}$ by
Proposition~\ref{prop3}. According to Corollary~\ref{cor2} (with\ \ \mb{E = \tilde{E}\ ,}\ \mb{D
= \tilde{D}\ ,}\ \mb{f = \id_E}\ )\ it holds\ \ \mb{k = l\ ,}\ and we may assume
that there are isomorphisms\ \ \mb{f_i : (E_i,D_i,h_i) \map
(\tilde{E}_i,\tilde{D}_i,\tilde{h}_i)}\ \ of flat bundles of $g$-degree 0 with
$g$-Einstein metrics, where\ \ \mb{D_i := D\vert_{E_i}\ ,}\
\mb{\tilde{D}_i := D\vert_{\tilde{E}_i}\ ,}\ \mb{h_i := h\vert_{E_i}\ ,}\
\mb{\tilde{h}_i := \tilde{h}\vert_{\tilde{E}_i}\ .}\ This means in particular
that the Higgs operator $d\dpr_i$ in $E_i$ associated to $D_i$ and $h_i$ is
isomorphic via $f_i$ to the Higgs operator $\tilde{d}\dpr_i$ in $\tilde{E}_i$
associated to $\tilde{D}_i$ and $\tilde{h}_i$. Hence\ \ \mb{d\dpr_h =
d\dpr_1\oplus\ldots d\dpr_k}\ \ is isomorphic to\ \ \mb{d\dpr_{\tilde{h}} =
\tilde{d}\dpr_1\oplus\ldots\oplus \tilde{d}\dpr_k}\ \ via the isomorphism\ \
\mb{f :=f_1\oplus\ldots\oplus f_k\ .}

In the same way, but using Proposition~\ref{prop10} and the results of section 4, one
shows that there is a well defined map from $\mathcal{M}\dpr(E)^0_g$ to
$\mathcal{M}_f(E)^0_g$, associating to the class of an integrable Higgs operator
$d\dpr$ with $g$-Einstein metric $h$ the class of the connection\ \ \mb{D_h =
I_h^{-1}(d\dpr)\ ;}\ this obviously is an inverse of $I$. \qed

\section{Line bundles on non-K\"ahler surfaces.}

Isomorphism classes of flat complex \underbar{line} bundles $(L,D)$ on a
manifold $X$ are parametrized by $H^1(X,\mbb{C}^*)$. On the other hand, an
integrable Higgs operator\ \ \mb{d\dpr = \dbar + \theta}\ \ in a complex
line bundle $L$ consists of a holomorphic structure $\dbar$ in $L$ and a
holomorphic 1-form $\theta$ on $X$ (the condition\ \
\mb{\theta\wedge\theta = 0}\ \ now is trivial). Furthermore, two integrable
Higgs operators $d\dpr_1$ and $d\dpr_2$ in $L$ are isomorphic if and only if
the two holomorphic line bundles $(L,\dbar_1)$ and $(L,\dbar_2)$ are isomorphic
and\ \ \mb{\theta_1 = \theta_2\ .}\ Hence, the space parametrizing isomorphism
classes of integrable Higgs operators is\ \ \mb{H^1(X,\mathcal{O}^*)\oplus
H^0(X,\Omega^1(X)) = \Pic(X)\oplus H^{1,0}(X)\ .} In particular, the moduli
sets $\mathcal{M}_f(L)^0_g$ and $\mathcal{M}\dpr(L)^0_g$ defined in the previous section are
subsets of $H^1(X,\mbb{C}^*)$ resp. $\Pic(X)\oplus H^{1,0}(X)$.

\begin{lem}\label{lem6} Let $L$ be a complex line bundle on $X$, and $g$ a Hermitian metric in $X$. Then every flat
connection in $L$ and every integrable
Higgs operator in $L$ admits a $g$-Einstein metric. \end{lem}

{\bf Proof:} Let $h_0$ be fixed metric in $L$, then every metric is of the
form\ \ \mb{h_f = e^f\cdot h_0}\ \ with
\ \ \mb{f \in \mathcal{C}^\infty(X,\mbb{R})\ .} Let
$D$ be a flat connection in $L$; then $h_f$ is a $g$-Einstein 
metric for $D$ if and only if it is a solution of the 
equation\ \ \mb{i\Lambda_gG_{h_0} - {i\over
2}\Lambda_g\dbar\partial(f) = c}\ \ with a real constant $c$. 
Such a solution exists by [LT] Corollary 7.2.9.
A similar argument works for integrable Higgs operators. \qed

{}From now on let $X$ be a \underbar{surface}, and $g$ a fixed Hermitian metric
in $X$. Then the map\ \ \mb{\deg_g : \Pic(X) \map \mbb{R}}\ \ is a morphism of Lie
groups ([LT] Proposition 1.3.7; recall that\ \ \mb{\deg_g = \deg_{\tilde{g}}}\ \ for some Gauduchon metric $\tilde{g}$). We define
\[
H^1(X,\mbb{C}^*)^f := \sett{[(L,D)] \in H^1(X,\mbb{C}^*)}{\deg_g(D) = 0}\ ,
\]
\[
\Pic(X)^T := \sett{[(L,\dbar)] \in \Pic(X)}{c_1(L)_\mbb{R} = 0}\ ,
\]
and 
\[
\Pic(X)^f := \ker(\deg_g\vert_{\Pic(X)^T})\ . 
\]
Observe that $\Pic(X)^f$ can be identified with the set of isomorphism
classes of line bundles admitting a flat \underbar{unitary} connection ([LT]
Proposition 1.3.13).

Theorem~\ref{thm3} and Lemma~\ref{lem6} imply

\begin{prop}\label{prop11} There is a natural bijection
\[
I_1 : H^1(X,\mbb{C}^*)^f \map \Pic(X)^f\times H^{1,0}(X)\ .
\]
\end{prop}

If $X$ admits a K\"ahler metric, i.e. if the first Betti number of $b_1(X)$ is
even, then $\deg_g$ is a topological invariant for every  metric
$g$ ([LT] Corollary 1.3.12 i)). Hence in this case it holds\ \ \mb{H^1(X,\mbb{C}^*)^f
= H^1(X,\mbb{C}^*)}\ \ and\ \ \mb{\Pic(X)^f = \Pic(X)^T\ ,}\ and $I_1$ is the
natural bijection from the moduli space of
isomorphism classes of flat line bundles to the moduli space of integrable
Higgs operators in line bundles with vanishing first real Chern class, which
(e.g. by the work of Simpson) already is known to exist for a K\"ahler metric $g$.

So let us assume that $b_1(X)$ is odd. Then\ \ \mb{\deg_g\vert_{\Pic^0(X)} :
\Pic^0(X) \map \mbb{R}}\ \ is surjective, and it holds
\[
\qmod{\Pic(X)^T}{\Pic(X)^f} \cong \qmod{\Pic^0(X)}{\Pic^0(X)^f} \cong \mbb{R}
\]
([LT] Corollary 1.3.12 and Proposition 1.3.13). We will show that $I_1$
extends to a (non-natural) bijection from $H^1(X,\mbb{C}^*)$ to $\Pic(X)^T\times
H^{1,0}(X)$ in this case, too. 

\begin{lem}\label{lem7} There is a bijection\ \ \mb{i : \Pic(X)^T \map
\Pic(X)^f\times\mbb{R}}\ \ such that the diagram
\[
\begin{array}{ccc}
\Pic(X)^T & \ovrightarrow{\deg_g} & \mbb{R} \\ 
\phantom{.}\\
i\downarrow\phantom{i} && \Vert \\
\phantom{.}\\
\Pic(X)^f\times\mbb{R} & \ovrightarrow{proj.} & \mbb{R}
\end{array}
\]
commutes.
\end{lem}

{\bf Proof:} $\deg_g\vert_{\Pic^0(X)}$ is surjective, so we can choose\ \ \mb{\mathcal{L}_1 := [(L_1,\dbar_1)] \in \Pic^0(X)}\ \ with
\ \ \mb{\deg_g(\mathcal{L}_1) = \deg_g(\dbar_1) = 1\ ,}\ and a class\ \ \mb{\alpha
\in H^1(X,\mathcal{O})}\ \ such that\ \ \mb{\mathcal{L}_1 = \pi(\alpha)}\ \ where\ \
\mb{\pi : H^1(X,\mathcal{O}) \map \Pic^0(X)}\ \ is the natural surjection. For\ \
\mb{\lambda \in \mbb{R}}\ \ define
\[
\mathcal{L}_\lambda :=  \pi(\lambda\cdot \alpha)\ ;
\]
then\ \ \mb{\deg_g(\mathcal{L}_\lambda) = \lambda}\ \
since\ \ \mb{\deg_g\circ\pi : H^1(X,\mathcal{O}) \map \mbb{R}}\ \ is linear. Now
define $i$ by
\[
i(\mathcal{L}) := (\mathcal{L}\otimes\mathcal{L}_{-\deg_g(\mathcal{L})},\deg_g(\mathcal{L}))\ ;
\]
then it is obvious that the inverse of $i$ is given by\ \ \mb{(\mathcal{L},\lambda)
\mapsto \mathcal{L}\otimes\mathcal{L}_\lambda\ ,}\ and that the diagram above commutes. \qed

In the proof of a similar statement for $H^1(X,\mbb{C}^*)$ we will use

\begin{lem}\label{lem8} The natural map       
\[
l^1 : H^1(X,\mbb{C}^*) \map \Pic(X)^T\ \ ,\ \ l^1([(L,D)]) := [(L,D\dpr)]\ .
\]
is surjective, i.e. a holomorphic structure $\dbar$ in a differentiable line
bundle $L$ on $X$ is the (0,1)-part of a flat connection if and only if the
real first Chern class $c_1(L)_\mbb{R}$ vanishes. \end{lem}

{\bf Proof:} $\Pic(X)^f$ can be naturally identified with $H^1(X,U(1))$, such
that the inclusion\ \ \mb{\Pic(X)^f \hookrightarrow \Pic(X)}\ \ becomes the
injection\ \ \mb{k^1 : H^1(X,U(1)) \hookrightarrow H^1(X,\mathcal{O}^*)}\ \ ([LT] p. 38).
$k^1$ is the composition of the natural map\ \ \mb{H^1(X,U(1)) \map H^1(X,\mbb{C}^*)}\ \ and
$l^1$, so it holds 
\[
\Pic(X)^f =\im(k^1) \subset \im(l^1)\ .
\]
Each component of $\Pic(X)^T$ contains a component of $\Pic(X)^f$ ([LT]
Remark 1.3.10), hence for each component
\[
\Pic^c(X) := \sett{[(L,\dbar)] \in \Pic(X)}{c_1(L)_\mbb{Z} = c} \subset \Pic(X)^T
\]
there exists a class\ \ \mb{[(L_c,D_c)] \in
H^1(X,\mbb{C}^*)}\ \ such that\ \ \mb{l^1([(L_c,D_c)]) \in \Pic^c(X)\ .}\ Define\
\ \mb{H^1(X,\mbb{C}^*)^0 := \sett{[(L,D)] \in H^1(X,\mbb{C}^*)}{c_1(L)_\mbb{Z} = 0}\ .}\
The commutative diagram with exact rows
\[
\begin{array}{ccccccccc}
0&\map &\mbb{Z} &\map &\mbb{C} & \ovrightarrow{exp}&\mbb{C}^* &\map &0\\
\phantom{.}\\
&&\Vert &&\downarrow&&\downarrow \\
\phantom{.}\\
0&\map &\mbb{Z} &\map &\mathcal{O} & \ovrightarrow{exp}&\mathcal{O}^* &\map &0
\end{array}
\]
induces the commutative diagram 
\[
\begin{array}{ccc}
H^1(X,\mbb{C}) &\map&H^1(X,\mbb{C}^*)^0\\
\phantom{.}\\
h^1 \downarrow\phantom{h^1} && \phantom{l^1}\downarrow l^1\\
\phantom{.}\\
H^1(X,\mathcal{O}) &\map& \Pic^0(X)
\end{array}
\]
with surjective horizontal arrows. Since $X$ is a surface, the left vertical
arrow $h^1$ is also surjective ([BPV] p. 117), hence $l^1$ maps $H^1(X,\mbb{C}^*)^0$
surjectively onto $\Pic^0(X)$. Now it is easy to see that every element of\ \
\mb{\Pic^c(X) \subset \Pic(X)^T}\ \ is of the form $l^1([(L_c\otimes
L,D_c\otimes D)])$ for some\ \ \mb{[(L,D)] \in H^1(X,\mbb{C}^*)^0\ .} \qed

\begin{lem}\label{lem9} There is a bijection\ \ \mb{j : H^1(X,\mbb{C}^*) \map
H^1(X,\mbb{C}^*)^f\times\mbb{R}}\ \ such that the diagram
\[
\begin{array}{ccc}
H^1(X,\mbb{C}^*) & \ovrightarrow{\deg_g^\prime} & \mbb{R} \\ 
\phantom{.}\\
j\downarrow\phantom{j} && \Vert \\
\phantom{.}\\
H^1(X,\mbb{C}^*)^f\times\mbb{R} & \ovrightarrow{proj.} & \mbb{R}
\end{array}
\]
commutes, where\ \ \mb{\deg_g^\prime := \deg_g\circ l^1\ }\ is the map
associated to the $g$-degree of flat connections. \end{lem}

{\bf Proof:} Choose\ \ \mb{\mathcal{L}_1 \in \Pic^0(X)\ ,}\ \mb{\alpha \in
H^1(X,\mathcal{O})}\ \ as in the proof of Lemma~\ref{lem7}, and a class\ \ \mb{\beta \in
H^1(X,\mbb{C})}\ \ with\ \ \mb{h^1(\beta) = \alpha\ .}\ Let\ \
\mb{\pi^\prime : H^1(X,\mbb{C}) \map H^1(X,\mbb{C}^*)}\ \ be the map induced by\ \
\mb{exp : \mbb{C} \map \mbb{C}^*\ ,}\ and define\ \
\mb{\mathcal{L}_1^\prime := \pi^\prime(\beta) \in H^1(X,\mbb{C}^*)\ .}\ Since the diagram
\[
\begin{array}{ccc}
H^1(X,\mbb{C}) & \ovrightarrow{\pi^\prime} & H^1(X,\mbb{C^*}) \\ 
\phantom{.}\\
h^1\downarrow\phantom{h^1} && \phantom{l^1}\downarrow l^1\\
\phantom{.}\\
H^1(X,\mathcal{O}) & \map & \Pic(X)^T
\end{array}
\]
commutes, it holds\ \ \mb{\deg_g^\prime(\mathcal{L}_1^\prime) = 1\ .}\ The rest of
the proof is as for Lemma~\ref{lem7}. \qed 

We conclude

\begin{thm} The composition
%
\[
\bar{I} : H^1(X,\mbb{C}^*) \ovrightarrow{j} H^1(X,\mbb{C}^*)^f\times\mbb{R}
\ovrightarrow{I_1\times \id_\mbb{R}} H^{1,0}(X)\times \Pic(X)^f\times\mbb{R} \]
\[\ovrightarrow{\id_{H^{1,0}(X)}\times i^{-1}} H^{1,0}(X)\times\Pic(X)^T\]
%
is a bijective extension of the map $I_1$, and preserves the $g$-degree.
\end{thm}

We finish with the obvious remark that the map\ \ \mb{l^1 : H^1(X,\mbb{C}^*) \map
\Pic(X)^T}\ \ in general does not coincide with the composition of $\bar{I}$
and projection onto $\Pic(X)^T$.

\Addresses
\end{document}